\documentclass[12pt]{amsart}

\usepackage{graphicx} 
\usepackage{tabularx}
\usepackage{latexsym,amssymb,amsmath,amsfonts,amscd} 
\usepackage[all]{xy}
\usepackage{mathrsfs} 
\usepackage{tikz} 
\usetikzlibrary{shapes,patterns,arrows,decorations.pathreplacing}
\usepackage{epigraph,enumerate}
\usepackage{stmaryrd}
\usepackage{tikz-cd}
\usepackage{mathtools}
\usepackage{hyperref}

\usepackage{bbm,cancel}
\allowdisplaybreaks[4]
\usepackage{latexsym, amssymb, amsmath}
\usepackage{enumerate}
\usepackage{flushend} 
\usepackage{mathrsfs} 
\usepackage{stfloats}
\def\custombibliography#1{
 \normalsize
\section*{\centering References}
 \list
 {[\arabic{enumi}]}{\settowidth\labelwidth{[#1]}\leftmargin\labelwidth
 \setlength{\itemsep}{.1em}
 \advance\leftmargin\labelsep
 \usecounter{enumi}}
 \def\newblock{\hskip .11em plus .33em minus -.07em}
 \sloppy
 \sfcode`\.=1000\relax}

\def\L2{{\cal L}_2}

\newcommand\bull{\vrule height .9ex width .8ex depth -.1ex } 

\newcommand\re{\rm I\! R}
\newcommand\cdcout[1]{} 

\newcommand{\rv}[1]{\boldsymbol{#1}} 
\newcommand{\RomanNumber}[1]{\uppercase\expandafter{\romannumeral #1}}
\newcommand{\romannumber}[1]{\lowercase\expandafter{\romannumeral #1}}

\DeclareMathAlphabet{\mathpzc}{OT1}{pzc}{m}{it}

\def\1{\rv 1} 

\usepackage{flushend} 
\usepackage{stfloats}
\usepackage{color} 
\usepackage{latexsym,amssymb,amsmath} 
\usepackage{mathrsfs,mathtools} 
\usepackage{multicol,enumerate}
\usepackage{subcaption}

\def\abs#1{\lvert #1 \rvert}

\def\allpoly{\mbox{$\re\langle X \rangle$}}

\def\allproperpoly#1{\mbox{$\re^{#1}_{p}\langle X \rangle$}}
\def\allpolyXpri{\mbox{$\re\langle X^{\prime} \rangle$}}

\def\allpolyx0degn{\mbox{$P_n$}}

\def\allwords{\mbox{$X^{\ast}$}}

\def\allseries{\mbox{$\re\langle\langle X \rangle\rangle$}}
\def\allseries#1{\mbox{$\re^{#1}\langle\langle X \rangle\rangle$}}
\def\allseriesXpri#1{\mbox{$\re^{#1}\langle\langle X^{\prime} \rangle\rangle$}}

\def\allseriesell{\mbox{$\re^{\ell} \langle\langle X \rangle\rangle$}}

\def\allseriesLC#1{\mbox{$\re^{#1}_{LC}\langle\langle X \rangle\rangle$}}

\def\allproperseries#1{\mbox{$\re_{p}^{#1}\, \langle\langle X \rangle\rangle$}}
\def\allnproperseries#1{\mbox{$\re_{np}^{#1}\, \langle\langle X \rangle\rangle$}}
\def\allpiseries#1{\mbox{$\re_{pi}^{#1}\, \langle\langle X \rangle\rangle$}}
\def\allpiseriesXpri#1{\mbox{$\re_{pi}^{#1}\, \langle\langle X^{\prime} \rangle\rangle$}}

\def\allseriesellLC{\mbox{$\re^{\ell}_{LC}\langle\langle X \rangle\rangle$}}

\def\allseriesX1{\mbox{$\re [[ X_1 ]]$}}

\def\bull{\rule{0.08in}{0.08in}} 

\newcommand{\comment}[1]{} 

\def\doubleone{{\rm\, l\!l}}

\def\End#1{{\rm End}\left(#1\right)}

\def\Endallseries{{\rm End}\left(\allseries{}\right)}

\def\eqref#1{(\ref{#1})} 

\def\lhook{\curvearrowleft}

\def\id{{\mbf{\rm id}}}
\def\ker{{\rm ker}}

\def\mbf#1{\hbox{\mathversion{bold}$#1$}} 
\def\nat{{\mathbb N}} 
\def\norm#1{\Vert#1\Vert}

\def\openbull{\framebox[0.08in][c]{$\;$}} 

\def\re{{\mathbb R}} 

\def\shuffle{{\scriptscriptstyle \;\sqcup \hspace*{-0.05cm}\sqcup\;}}

\def\supp{{\rm supp}}

\def\begals{\[\begin{aligned}}
\def\endals{\end{aligned}\]}
\def\begal{\begin{align*}}
\def\endal{\end{align*}}
\def\begce{\begin{center}}
\def\endce{\end{center}}
\def\begar{\begin{array}}
\def\endar{\end{array}}
\def\begeq{\begin{equation}}
\def\endeq{\end{equation}}
\def\begdi{\begin{displaymath}}
\def\enddi{\end{displaymath}}
\def\begdis{\begin{eqnarray*}}
\def\enddis{\end{eqnarray*}}
\def\begeqa{\begin{eqnarray}}
\def\endeqa{\end{eqnarray}}
\def\begdes{\begin{description}}
\def\enddes{\end{description}}
\def\begit{\begin{itemize}}
\def\endit{\end{itemize}}
\def\begen{\begin{enumerate}}
\def\enden{\end{enumerate}}
\def\beglar{\left[\begin{array}}
\def\endrar{\end{array}\right]}
\def\begle{\begin{lemma}}
\def\endle{\end{lemma}}
\def\begde{\begin{definition}}
\def\endde{\end{definition}}
\def\begth{\begin{theorem}}
\def\endth{\end{theorem}}
\def\begco{\begin{corollary}}
\def\endco{\end{corollary}}
\def\begprop{\begin{proposition}}
\def\endprop{\end{proposition}}	
\def\begex{\begin{example}}
\def\endex{\hfill\openbull \end{example}}
\def\begexer{\begin{exercise}}
\def\endexer{\end{exercise}}
\def\begalg{\begin{algo}}
\def\endalg{\end{algo}}
\def\begre{\noindent{\bf Remark}:\hspace*{0.05cm}}

\def\begres{\noindent{\bf Remarks}:\begin{enumerate}}
\def\endres{\end{enumerate} \par}
\def\begpr{\noindent{\em Proof:}$\;\;$}
\def\endpr{\hfill\bull}
\def\begtab{\begin{tabular}}
\def\endtab{\end{tabular}}
\def\rref#1{(\ref{#1})}
\newtheorem{lemma}{Lemma}[section]
\newtheorem{definition}{Definition}[section]
\newtheorem{theorem}{Theorem}[section]
\newtheorem{proposition}{Proposition}[section]
\newtheorem{corollary}{Corollary}[section]
\newtheorem{example}{Example}[section]
\newtheorem{algo}{Algorithm}[section]
\def\HAprod{\mbf{m}} 
\def\shuff#1#2{\mathbin{
		\hbox{\vbox{\hbox{\vrule \hskip#2 \vrule height#1 width 0pt}\hrule}\vbox{\hbox{\vrule \hskip#2 \vrule height#1 width 0pt\vrule }\hrule}}}}
\def\shuffl{{\mathchoice{\shuff{5pt}{3.5pt}}{\shuff{5pt}{3.5pt}}{\shuff{3pt}{2.6pt}}{\shuff{3pt}{2.6pt}}}}
\def\shuffle{{\, \shuffl \,}}
\def\shuff#1#2{\mathbin{
      \hbox{\vbox{\hbox{\vrule \hskip#2 \vrule height#1 width 0pt}\hrule}\vbox{\hbox{\vrule \hskip#2 \vrule height#1 width 0pt\vrule }\hrule}}}}
\def\shuffl{{\mathchoice{\shuff{5pt}{3.5pt}}{\shuff{5pt}{3.5pt}}{\shuff{3pt}{2.6pt}}{\shuff{3pt}{2.6pt}}}}
\def\shuffle{{\, \shuffl \,}}

\newcommand{\rvline}{\hspace*{-\arraycolsep}\vline\hspace*{-\arraycolsep}}

\begin{document}

\title[A Formal Series Approach to Multiplicative Dynamic Feedback Connection]{A Formal Power Series Approach to Multiplicative Dynamic Feedback Interconnection}

\author{Kurusch Ebrahimi-Fard}
\address{Department of Mathematical Sciences, Norwegian University of Science and Technology (NTNU), 7491 Trondheim, Norway}
\email{kurusch.ebrahimi-fard@ntnu.no}
\urladdr{https://folk.ntnu.no/kurusche/}

\author{Venkatesh~G.~S.}
\address{Department of Mathematical Sciences, Norwegian University of Science and Technology (NTNU), 7491 Trondheim, Norway}
\email{subbarao.v.guggilam@ntnu.no}


\begin{abstract} 
The goal of the paper is multi-fold. First, an explicit formula is derived to compute the non-commutative generating series of a closed-loop system when a (multi-input, multi-output) plant, given in Chen--Fliess series description is in multiplicative output feedback interconnection with another system, also given as Chen--Fliess series. Furthermore, it is shown that the multiplicative dynamic output feedback connection has a natural interpretation as a transformation group acting on the plant. A computational framework for computing the generating series for multiplicative dynamic output feedback is devised utilizing the Hopf algebras of the coordinate functions corresponding to the shuffle group and the multiplicative feedback group. The pre--Lie algebra in multiplicative feedback is shown to be an example of Foissy's com-pre-Lie algebras indexed by matrices with certain structure. \\

\noindent {\tiny{{\bf{Key words}}: Chen--Fliess series; feedback; Hopf algebras; pre- and post-Lie algebras; com-pre-Lie algebra.}}

\smallskip 

\noindent {\tiny{\textbf{MSC2020}: 
93C10 (primary),  	
16T05 (primary), 	
16T30 (secondary), 	
17D25 (secondary)}} 	
\end{abstract}

\maketitle

\noindent
\tableofcontents

\thispagestyle{empty}


\section{Introduction}
\label{sec:intro}

The document studies multiplicative feedback interconnections of input-output systems using the framework of Chen--Fliess functional series~\cite{Fliess_81}. As such, there is no need for a state space realization and thus, the results presented here are independent of any state space embedding when a realization is possible~\cite{Fliess_realizn_83}. It was shown in~\cite{Gray-Li_05} that {\em additive feedback} interconnection of two such systems results in a Chen--Fliess series description for the closed--loop system. An efficient computation of the generating series for such closed--loop systems is facilitated through use of tools from combinatorial Hopf algebra~\cite{Duffaut-Espinosa-etal_JA16,Gray-etal_SCL14}. 

\smallskip

{\it{Multiplicative feedback}} is another vital interconnection of dynamical systems that is extensively studied in control theory and communications~\cite{Haining_Cass_14,Isidori_95}. However, when the nature of feedback interconnection becomes {\em multiplicative}, a description in terms of combinatorial Hopf algebra has been unexplored. It is known that, in the single-input single-output (SISO) setting, the closed-loop system in the affine feedback case (of which multiplicative feedback is a special case) has a Chen--Fliess series description and the computation of feedback formula is facilitated through combinatorial Hopf algebra of coordinate functions~\cite{Gray-KEF_SIAM2017} (see also \cite{Foissy_2016}).
 
\smallskip

The contribution of the article at hand are two--fold, combining {\em systems theory} and {\em algebra}. Regarding the former, 
\begen[(i)]

\item It is shown that in the multi-input multi-output (MIMO) setting the closed-loop system under multiplicative feedback has a Chen--Fliess series representation. An explicit expression of the closed-loop generating series is provided, which will be called {\em multiplicative dynamic feedback product}. 

\medskip

\item The invariance of class and relative degree of a dynamical system (plant) under the multiplicative feedback connection is shown. The invariance of relative degree of the plant is essential in the design of feedback linearization of the closed--loop system, a question which is deferred for the future work.

\medskip

\item It will be shown that the closed--loop multiplicative feedback product has a natural interpretation as a transformation group acting on the plant. The transformation group, denoted by $\left(\allpiseries{m}, \shuffle\right)$, is the shuffle group on vectors with formal power series over non-commuting set of indeterminates $X$, whose constant term is $1$ in each coordinate, as entries. 

\medskip

\item Prior to discussing transformation group, it is essential to describe the unit group $\left(\allpiseries{m}, \star\right)$, that arises from cascade interconnection of Chen--Fliess series along with their multiplicative feedforward of inputs. It turns out to be the Grossman--Larson group of the pre-group $\left(\allpiseries{m}, \shuffle, \lhook\right)$. The details of action map $\lhook$ are to be found in Section~\ref{subsec:dyn_mult_feedb_grp}. For more details on pre-groups in general, we refer to~\cite{Alkaabi2023postgroup,Bai2022postgroup}.    
\enden

\smallskip

Part of the work described so far has appeared in earlier work of the first author \cite{Venkat_MTNS_22}. 

\smallskip

\noindent The contributions towards the algebra side are briefly tip-toed as follows:

\begen[(i)]
\item The algorithmic framework for the computation of the multiplicative dynamic feedback product formula for a general MIMO case is devised using a connected graded bialgebra $\mathcal{H}$, which is shown to be right-handed (in the sense of \cite{MenousPatras2015}).

\medskip

\item The character group of the Hopf algebra $\mathcal{H}$ is isomorphic to a normal subgroup of $\left(\allpiseries{m}, \star\right)$~(details in Section~\ref{sec:computational_framework}). 

\medskip

\item A cointeracting Hopf algebra structure $\left(\mathcal{H},\Delta_{\shuffle}, \Delta, \epsilon\right)$ is shown to result from dualizing the pre-group structure mentioned earlier. The coproducts $\Delta$ dualizes dynamic feedback group product denoted $\star$. The action map $\lhook$ dualizes to the coaction map $\rho$.

\medskip 

\item Linearizing the right-handed coproduct $\Delta$ results in a graded right pre-Lie coalgebra and its graded dual defines a graded right pre-Lie product on the space of non-commutative polynomials with constant term being zero in each coordinate, denoted by $\allproperpoly{m}$. Linearizing the coaction map and taking its graded dual gives a right pre-Lie product. In fact, it is shown that $\left(\allproperpoly{m}, \triangleleft\right)$ is a special case of Foissy's {\em com-pre-Lie} algebras.

\medskip

\item We show that over a field of characteristic zero, every com-pre-Lie algebra inherits a new  pre-Lie product (Theorem~\ref{th:com-pre-Lie-2-pre-Lie}). The two pre-Lie products anti-symmetrize to same Lie bracket; the multiplicative feedback case is a prime example of this observation.     
\enden

\medskip
\medskip

The paper is organized as follows. The next section provides a brief summary of Chen--Fliess series and their interconnections. Section~\ref{subsec:dyn_mult_feedb_grp} builds the pivotal {\em multiplicative dynamic output feedback group} that acts as a transformation group for the multiplicative feedback interconnection. Section~\ref{sec:mult_dyn_feedb} is where the closed--loop configuration of Chen--Fliess series is analyzed. The invariance of relative degree under multiplicative output feedback is asserted in Section~\ref{sec:rel_degree_mult_feedb}. A framework for computing the feedback product is devised in Section~\ref{sec:computational_framework} including examples. The discussion on pre-Lie products and com-pre-Lie algebras are carried out in Sections~\ref{sec:pre-Lie} \& \ref{sec:com-pre-Lie}

  
\subsection{Preliminaries: Formal Power Series}
\label{ssec:prelims}

We set notation by recalling some preliminaries about non-commutative formal power series and various algebraic structures, including bialgebras. For more details regarding the latter, we refer to~\cite{Abe_Hopf,CatPat2021,Manchon_08,Reutenauer_93,Sweedler}. 

\smallskip

\noindent Let $X = \{x_0,x_1,\ldots,x_m\}$ denote a finite alphabet of non-commuting indeterminates.

\begen

\item $\allpoly$ is the $\re$-algebra of non-commutative polynomials with indeterminates in $X$. As a vector space, it is spanned by the set of words, denoted by $\allwords$, with letters from $X$. The unit is the empty word, which we write $\mathbf{1}$, and the set of non-empty words is denoted by $X^+:=X^*\backslash\{\mathbf{1}\}$.
 \medskip
 
\item $\allpoly$ is a cocommutative Hopf algebra, graded (by the word length, $|w| \ge 0$, $w \in X^*$) and connected, with catenation product and counit $\epsilon$ such that $\epsilon(w) = \delta_{|w|,0}$ for any $w \in X^*$. The coproduct is unshuffle $\Delta_{\shuffle}$ for which the elements $x_i \in X$ are primitive, i.e., $\Delta_{\shuffle}(x_i) = x_i \otimes \mathbf{1} + \mathbf{1} \otimes x_i$. A few examples of this coproduct are given next
\begin{align*}
	\Delta_{\shuffle}x_ix_j 
	&= x_ix_j \otimes \mathbf{1} + x_i \otimes x_j + x_j \otimes x_i + \mathbf{1} \otimes x_ix_j, \\
	\Delta_{\shuffle}x_i^2 
	&= x_i^2 \otimes \mathbf{1} + 2 x_i \otimes x_i + \mathbf{1} \otimes x_i^2. 
\end{align*}
In general, with $ x_i^0 :=  \mathbf{1}$, we have
\begin{align}
\label{eqn:shuffle_pow}
	\Delta_{\shuffle} x_i^n 
	&= \sum\limits_{k = 0}^n \binom{n}{k} x_i^k \otimes x_i^{n-k}.
\end{align}

\item \medskip

\begre{
	Traditionally, the counit of $\allpoly$ is denoted by $\emptyset$ in the Chen--Fliess literature; for the rest of the article the presentation also follows this custom.
}

\medskip

\item The convolution algebra of linear maps from the counital coalgebra $\left(\allpoly, \Delta_{\shuffle}\right)$ to $\re$, is given by the space of formal power series denoted by $\allseries{}$. The dual basis is given by $\{\emptyset\} \cup X^{+}$, such that $\eta\left(\xi\right) = 1$ if $\eta = \xi$ in $X^+$, and zero else. An element $c \in \allseries{}$ is represented by
\begin{align*}
	c = c(\mathbf{1})\emptyset + \sum\limits_{\eta \in X^+} c(\eta) \eta.
\end{align*} 
In the following, $\emptyset$ is not explicitly written unless needed. 
\medskip

The convolution product on $\allseries{}$ is the shuffle product, which is defined for all $c,d \in \allseries{}$ and $p \in \allpoly$ by
\begin{align*}
	\left(c \shuffle d\right)(p) 
	&= m_{\re} \circ \left(c \otimes d\right) \circ \Delta_{\shuffle} (p).
\end{align*} 
Here, $m_{\re}$ is the usual product in $\re$ and the unit element is the counit $\emptyset$. 
\medskip

\item The shuffle group $\left(M, \shuffle\right)$ is the normal subgroup of the unit group of $\allseries{}$ where 
\begin{align*}
	M &= \{1 \emptyset + c': c' \in \allseries{},  c'( \mathbf{1}) = 0\}.
\end{align*}

\item Consider arbitrary $c,d \in \allseries{}$ and $p,q \in \allpoly$. The (left) shift operation $p^{-1}$ on the shuffle algebra $\left(\allseries{},\shuffle\right)$ is defined by $p^{-1}(c)(q) := c(pq)$. For primitive elements $x_i \in X \subsetneq \allpoly$, the left shift is a derivation on the shuffle algebra viz.
\begin{align*}
x_i^{-1}\left(c \shuffle d\right) &=  \left(x_{i}^{-1}\left(c\right)\right) \shuffle d + c \shuffle \left(x_i^{-1}\left(d\right)\right). 
\end{align*}
\enden

\medskip

Looking ahead, we define for $c \in \allseries{}$ the Chen--Fliess series $F_c[u](t)$ as a functional series such that the basis (words) of $\allseries{}$ are identified with iterated integrals of the function $u$, a control signal (detailed in Section~\ref{sec:Interconnections}). In the context of a dynamical system, a Chen--Fliess series expresses its input-output behavior and provides a coordinate-independent framework to analyze the intrinsic properties of the system.

\label{ssec:facts}

\begen
\item $\allseries{}$ can also be considered a non-commutative $\re$-algebra under catenation product, defined for all $c,d \in \allseries{}$ and $\eta \in \allwords$ by
\begin{align*}
	c.d\left(\eta\right) = \sum\limits_{\substack{\zeta,\mu \in X^* \\ \zeta\mu = \eta}} c\left(\zeta\right)d\left(\mu\right).
\end{align*}
\medskip

\item An element $c \in \allproperseries{} := \{c \in \allseries{} : c(\mathbf{1}) = 0\}$ is termed {\em proper} series. Whereas $\allnproperseries{} := \allseries{}\setminus\allproperseries{}$ contains non-proper series. An element $c \in \allproperseries{}$ which has {\em finite support}, is termed {\em proper} polynomial. The set of proper polynomials are denoted by $\allproperpoly{}$. 

\medskip

\item $\allseries{}$ is a complete topological vector space endowed by ultrametric $\kappa$ defined as $\kappa (c,d) = \sigma^{-val\left(c-d\right)},$ for all $c,d \in \allseries{}$ and $\sigma \in ]0,1[$. Here $val\left(\cdot\right)$ is the valuation of a series.  See \cite{Berstel-Reutenauer_88} for details.

\medskip
	
\item Consider the $\ell$-fold Cartesian product $\allseries{\ell} \; \cong \; \allseries{} \times \allseries{} \times \cdots \times \allseries{}$ as complete topological vector space. Thus, a Chen--Fliess series $F_c$ for $c = \left(c_1,c_2,\ldots,c_{\ell}\right)\in \allseries{\ell}$ is a vector of Chen--Fliess series $\left(F_{c_1},F_{c_2},\ldots,F_{c_{\ell}}\right)$.

\medskip

\item $\left(\allseries{\ell}, \shuffle\right)$ is a direct product of shuffle algebras, $ \left(c \shuffle d\right)_j = c_j \shuffle d_j$, $c,d \in \allseries{\ell}$. 

\medskip

\item Following Foissy \cite{Foissy_2016}, $\allseries{\ell}$ can also inherit associative (but non-commutative) algebra structure via ``adorned" shuffle products, denoted $\shuffle_k$, where the subscript $k = 1,2,\ldots, \ell$. For $c,d \in \allseries{\ell}$
\begin{align*}
c \shuffle_k  d  =  \begin{pmatrix}
c_1 \shuffle d_k \\ c_2 \shuffle d_k \\ \vdots \\ c_{\ell} \shuffle d_k
\end{pmatrix}.
\end{align*}
\medskip

\item In general, for a given $\mathsf{a} = \left(a_1,a_2,\ldots,a_\ell\right) \in \re^{\ell}$; define the weighted product $c \shuffle_{\mathsf{a}} d = \sum_{i=1}^{\ell} a_i \left(c \shuffle_i d\right)$, then $\left(\allseries{\ell}, \shuffle_{\mathsf{a}}\right)$ is an associative (but non-commutative) algebra~\cite{Foissy_2016}. It follows from the next identity
\begin{align*}
c \shuffle_i \left(d \shuffle_j e\right) = \left(c \shuffle_i d\right) \shuffle_j e = \left(c \shuffle_j e\right) \shuffle_i d
\qquad c,d,e \in \allseries{\ell}. 
\end{align*}
\medskip

\item The products $\shuffle$ and $\shuffle_i$ are {\em mixed} associative, a property that is particularly made use of in Section~\ref{sec:com-pre-Lie} viz for all $c,d,e \in \allseries{\ell}$:
\begin{align*}
\left(c \shuffle_i d\right) \shuffle e = \left(c \shuffle e \right) \shuffle_i d 
\quad \text{ and } \quad
\left(c \shuffle_i d\right) \shuffle_i e = c \shuffle_i \left(d \shuffle e\right).
\end{align*}
\medskip
  
\item Consider $\allproperseries{\ell} \; \cong \; \allproperseries{} \times \allproperseries{} \times \cdots \times \allproperseries{}$ as topological vector subspace. The following set notation is thus imminent:
\begin{align*}
	\allpiseries{\ell} 
	&=  \allnproperseries{} \times \allnproperseries{} \times 
					\cdots \times \allnproperseries{}
\end{align*} 
and $c \in \allpiseries{\ell} $ is called {\em purely improper}. As $\allnproperseries{\ell} = \allseries{\ell}\setminus\allproperseries{\ell}$,  purely improper is thus a stronger notion than just non--proper when $\ell > 1$.
\medskip

\item A series $c \in \allseries{\ell}$ can be decomposed as 
$$
	c = c_N + c_F,
$$ 
where $c_N \in \mathbb{R}^\ell\langle\langle X_0 \rangle\rangle$, $X_0:=\{x_0\}$, is the so-called {\em natural part} of $c$ and $c_F:=c-c_N $ is termed its {\em forced part}. 
\medskip

\item $\re^{\ell}$ is treated as a commutative algebra under pointwise product with unit element $\doubleone = (1,1,\ldots,1)$.
\enden


\section{Chen--Fliess Series and the Interconnections}\label{sec:Interconnections}
\subsection{Chen--Fliess Series}

Let $\mathfrak{p}\ge 1$ and $T > 0$ be given. For a Lebesgue measurable function $u: [0,T] \rightarrow\re^m$, define
$\norm{u}_{\mathfrak{p}}=\max\{\norm{u_i}_{\mathfrak{p}}: \ 1 \le i \le m\}$, where $\norm{u_i}_{\mathfrak{p}}$ is the usual
$L_{\mathfrak{p}}$-norm for a measurable real-valued function, $u_i$, defined on $[0,T]$.  Let $L^m_{\mathfrak{p}}[0,T]$
denote the set of all measurable functions defined on $[0,T]$ having a finite $\norm{\cdot}_{\mathfrak{p}}$ norm and
$B_{\mathfrak{p}}^m(R)[0,T]:=\{u \in L_{\mathfrak{p}}^m[0,T]:\norm{u}_{\mathfrak{p}}\leq R\}$. 

\medskip

Given any series $c\in\allseriesell$, the corresponding
{\em Chen--Fliess series} is
\begeq
	F_c[u](t) = \sum_{\eta\in X^{\ast}} c(\eta)\, F_\eta[u](t), \label{eq:Fliess-operator-defined}
\endeq
where $F_{\emptyset}[u]:=1$ and
\begdi
	F_{x_i\bar{\eta}}[u](t) 
	:=\int_{0}^tu_{i}(\tau)F_{\bar{\eta}}[u](\tau)\,d\tau
\enddi
with $x_i\in X$, $\bar{\eta}\in X^{+}\cup \{\emptyset\}$, and $u_0:=1$ \cite{Fliess_81}. If there exist constants $K,M>0$ such that
\begin{align}
	\abs{c_i(\eta)}\leq K M^{\abs{\eta}} \abs{\eta}!,\;\; 
	\label{eq:Gevrey-growth-condition}
\end{align}
for all $\eta\in X^\ast$ and $i = 1,\ldots,\ell$, then $F_c$ constitutes a well-defined mapping from $B_{\mathfrak p}^m(R)[0, T]$ into $B_{\mathfrak q}^{\ell}(S)[0, \, T]$ for sufficiently small $R,T >0$, where the numbers $\mathfrak{p},\mathfrak{q}\in[1,\infty]$ are conjugate exponents, i.e., $1/\mathfrak{p}+1/\mathfrak{q}=1$~\cite{Gray-Wang_SCL02}. A series $c \in \allseriesell$ with coefficients obeying the growth condition in ~\rref{eq:Gevrey-growth-condition} is called a {\em locally convergent} generating series. The set of all locally convergent generating series is denoted by $\allseriesellLC$. The details regarding convergence of Chen--Fliess series with respect to the growth of their generating series are well-detailed in the literature~\cite{GS_thesis,Irina_thesis,Winter_Arboleda-etal_2015}. The expansion \eqref{eq:Fliess-operator-defined} is called a {Fliess operator} if $c \in \allseriesLC{\ell}$ otherwise it defines an operator only in a formal sense.\\

\begre{~The shuffle product on $\allseries{}$ corresponds to encoding the {\em integration by parts} rule of the Lebesgue--Stieltjes integration~\cite{Ree_58}.}

\subsection{Interconnections of Chen--Fliess Series}
\label{subsec:interconnections}

Given two Chen--Fliess series, $F_c$ and $F_d$, $c,d\in\allseriesell$, the parallel and product interconnections \cite{Fliess_81,Ree_58}, $F_c+F_d=F_{c+d}$ and $F_cF_d=F_{c\shuffle d}$ respectively, both preserve local convergence \cite{Thitsa-Gray_SIAM12,GS_thesis}. For $c \in \allseriesXpri{k}$ and $d \in \allseries{\ell}$, where $X^{\prime}= \{x^{\prime}_0,x^{\prime}_1,\ldots,x^{\prime}_\ell\}$, the composition, or cascade interconnection, $F_c \circ F_d$, has a Chen--Fliess series representation $F_{c \circ d} \in \allseries{k}$, defined via the {\em composition product} \cite{Ferfera_79,Ferfera_80} where
\begeq 
\label{eq:c-circ-d}
	c \circ d= \sum_{\eta\in X^{\prime\ast}} c(\eta)\,\psi_d(\eta)(\emptyset).
\endeq
The continuous (in ultrametric topology) algebra homomorphism $\psi_{d} : \allpolyXpri \rightarrow \Endallseries$ maps $x_i^{\prime} \eta \longmapsto \psi_d(x_i^{\prime}\eta)$ and is defined for all $i=0,1,\ldots,\ell$ and any $e \in \allseries{}$ by
$$
	\psi_d(x_i^{\prime})(e)= x_0(d_i \shuffle e), \qquad d_0:= \emptyset.
$$
By definition, $\psi_d(\emptyset)$ is the identity map on $\allseries{}$. The cascade interconnection preserves local convergence~\cite{Thitsa-Gray_SIAM12} and the composition product \eqref{eq:c-circ-d} distributes over the shuffle product.

\begth\cite{Gray-Li_05}
\label{thm:shuff_dist_comp} 
Let $c,d \in \allseriesXpri{k}$ and $e \in \allseries{\ell}$, such that $\abs{X^{\prime}} = \ell + 1$, then 
\begin{equation}
\label{eq:distrib}
	\left(c\shuffle d\right)\circ e = \left(c \circ e\right) \shuffle \left(d \circ e\right).
\end{equation} 
\endth

\medskip

The composition product preserves the purely improper property of the left argument.

\begth
\label{thm:proper_preser_comp} 
If $c \in \allseriesXpri{k}$ and $d \in \allseries{\ell}$ such that $\abs{X^{\prime}} = \ell +1$, then $\left(c\circ d\right) \left(\mathbf{1}\right) = c\left(\mathbf{1}\right)$. Hence, $c \in \allpiseriesXpri{k}$ if and only if $c\circ d \in \allpiseries{k}$ and vice-versa. Similarly, if $c$ is a proper series then $c \circ d$ is a proper series and vice-versa. 
\endth

\begpr 
The proof follows immediately from \rref{eq:c-circ-d}.
\endpr

\medskip

The composition product \rref{eq:c-circ-d} is a strong contraction map with respect to its right argument in the ultrametric topology.

\begth\cite{Gray-Li_05}
\label{thm:comp_strong_cont} 
Let $c \in \allseriesXpri{k}$ and $d,e \in \allseries{\ell}$, such that $\abs{X^{\prime}} = \ell + 1$, then $\kappa\left(c\circ d, c\circ e\right) \leq \sigma \kappa\left(d,e\right)$ where $\sigma \in [0,1[$.
\endth


\subsection{Cascading with Multiplicative Feedforward of Input}
\label{subsec:Chen_fliess_interconn} 

\begin{figure}[tb]
	\begin{center}
		\includegraphics[scale=0.75]{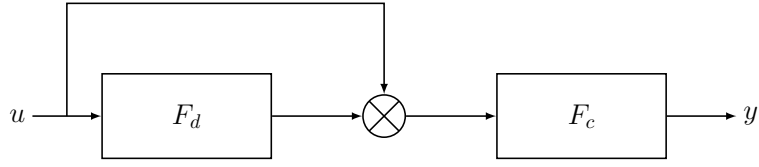}
		\caption{Cascade interconnection of Chen--Fliess series $F_d$ with $F_c$ along with multiplicative feedforward of input}
		\label{fig:mult_mix_comp}
	\end{center}
\end{figure}

The cascade interconnection of Chen--Fliess series along with the multiplicative feedforward of the input, as shown in Figure~\ref{fig:mult_mix_comp}, arises primarily in the analysis of multiplicative feedback interconnection discussed in Section~\ref{sec:mult_dyn_feedb}. A semblance of definition of such an interconnection has appeared implicitly (and was limited to the SISO case) in Definition~$3.1$ of ~\cite{Gray-KEF_SIAM2017}. 

\smallskip

The input--output map of Figure~\ref{fig:mult_mix_comp}, the function
$$
	y = F_c[u.F_d[u]]
$$ 
has Chen--Fliess series representation denoted by $F_{c \lhook d}$, where the {\em multiplicative mixed composition product} of $c \in \allseries{p}$ and $d \in \allseries{m}$ is defined by:
\begin{align}
\label{eqn:mixprocomp_def}
	c \lhook d 
	:= \sum_{\eta\in X^\ast} c\left(\eta\right) \eta \lhook d,
\end{align}
where for all $\eta \in \allwords$ and $d \in \allseries{m}$ the binary operation $\eta \lhook d$ is defined inductively
\begin{align*}
	\emptyset \lhook d 
	&:= \emptyset\\
	x_0\eta \lhook d 
	&:= x_0\left(\eta \lhook d\right)\\
	x_i\eta \lhook d 
	&:= x_i\left(d_i \shuffle (\eta \lhook d)\right) \quad \qquad \forall \, i = 1,2,\ldots,m.
\end{align*}

The following results were proven in the SISO setting~\cite{Gray-KEF_SIAM2017}. Their MIMO extensions are straightforward and thus attributed to the referenced document with the  omission of proofs.

\smallskip
 
The multiplicative mixed composition \eqref{eqn:mixprocomp_def} distributes over the shuffle product on the right.

\begth\cite{Gray-KEF_SIAM2017}
\label{thm:shuff_dist_mixprc} 
Let $c,d \in \allseries{p}$ and $e \in \allseries{m}$, then 
$$
	\left(c \shuffle d\right) \lhook e = \left(c \lhook e\right) \shuffle \left(d \lhook e\right).
$$
\endth
\medskip

It is a strong contraction in its right argument in the ultrametric topology.

\begth\cite{Gray-KEF_SIAM2017}
\label{thm:mixpr_cont} 
Let $d,e \in \allseries{m}$ and $c \in \allseries{p}$. Define $c^{\prime} = c - c(\mathbf{1})$, then 
$$
	\kappa\left(c \lhook d, c \lhook e\right) \leq \sigma^{val\left(c^{\prime}\right)}\kappa\left(d,e\right).
$$
\endth  

Observe that $val\left(c^{\prime}\right) \geq 1$ and $\sigma \in ]0,1[$ in Theorem~\ref{thm:mixpr_cont}. The composition product \eqref{eq:c-circ-d} and the multiplicative mixed composition product \eqref{eqn:mixprocomp_def} are associative in combination.

\begth\cite{Gray-KEF_SIAM2017}
\label{thm:mix_assoc_comp_mixpr} 
Let $X^{\prime} = \{x'_0,\ldots,x'_p\}$ and $c \in \allseriesXpri{q}$. Let $d \in \allseries{p}$ and $e \in \allseries{m}$, then 
$$
	c \circ \left(d \lhook e\right) = \left(c \circ d\right)\lhook e.
$$
\endth


\subsection{Multiplicative Dynamic Output Feedback Group}
\label{subsec:dyn_mult_feedb_grp}

The dynamic multiplicative feedback group plays a vital role in computation of the multiplicative dynamic feedback formula and appeared implicitly (in the simpler SISO setting) in reference \cite{Gray-KEF_SIAM2017}. 

\smallskip

\begin{figure}[tb]
	\begin{center}
		\includegraphics[scale=0.65]{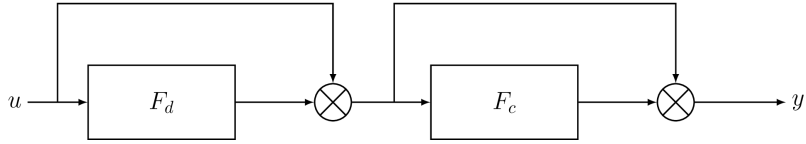}
		\caption{Cascade interconnection of Chen--Fliess series $F_d$ with $F_c$ along with multiplicative feedforward of their inputs.}
		\label{fig:mult-group}
	\end{center}
\end{figure}

Figure~\ref{fig:mult-group} shows the cascade interconnection of $F_c$ and $F_d$,  $c,d \in \allseries{m}$, along with their multiplicative feedforward of inputs. The input-output relation of the composite system is 
\begin{align}
\label{eqn:grp_prodF}
	y=u.F_d[u]F_c[u.F_d[u]] .
\end{align}
Defining the {\em multiplicative composition product} of $c,d \in \allseries{m}$
\begin{align}
\label{eqn:grp_prod}
	c \star d := d \shuffle \left(c \lhook d \right),
\end{align}
the right--hand side of \eqref{eqn:grp_prodF} can be represented by Chen--Fliess series, $y=u.F_{c \star d}[u].$ 

\begth\cite[Lem.~3.6]{Gray-KEF_SIAM2017}
\label{thm:mult_comp_assoc} 
Let $c,d,e \in \allseries{m}$, then $c \star \doubleone = c$ and
$$
	(c \star d )\star e = c \star (d \star e).
$$ 
\endth

We have a  right action by the non-commutative monoid $\left(\allseries{m},\star, \doubleone\right)$ on $\allseries{q}$.

\begth\cite{Gray-KEF_SIAM2017}
\label{thm:mixprocomp_right_act_monoid} 
Let $c \in \allseries{q}$ and $d,e \in \allseries{m}$, then 
\begin{align}
\label{eqn:rightact}
	\left(c \lhook d\right)\lhook e = c \lhook \left(d \star e\right).
\end{align}
\endth

The prominent question is to find the invertible elements of the monoid $\left(\allseries{m},\star, \doubleone\right)$ and the motivation to find the unit elements of the monoid shall be evident in Section~\ref{sec:mult_dyn_feedb}. Let $d,e \in \allpiseries{m}$ and suppose
$d \star e = \doubleone.$ Observe that $d \in \allpiseries{m}$ implies $\left(d \lhook e\right) \in \allpiseries{m}$ and using Theorem~\ref{thm:shuff_dist_mixprc} yields
\begin{align*}
	e = \left(d \lhook e \right)^{\shuffle -1} = d^{\shuffle -1} \lhook e.
\end{align*}
Thus for $e$ to be right inverse of $d$, it has to satisfy the fixed point equation 
\begin{align}
\label{eqn:fixed_point_grp_eqn_rinv}
	e = d^{\shuffle -1} \lhook e.
\end{align}
Theorem~\ref{thm:mixpr_cont} asserts that the map $e \mapsto d^{\shuffle -1} \lhook e$ is a strong contraction in the ultrametric space inferring that \rref{eqn:fixed_point_grp_eqn_rinv} has a unique fixed point. 
\medskip

Suppose $e$ is the left inverse of $d$, i.e., $e \star d = \doubleone$, then it satisfies the equation,
\begin{align}
\label{eqn:grp_eqn_linv}
	d = e^{\shuffle -1} \lhook d.
\end{align}
If $e$ is a solution of \rref{eqn:fixed_point_grp_eqn_rinv}, then it satisfies \rref{eqn:grp_eqn_linv} and vice--versa. Thus, $e$ is the unique inverse of $d \in \allpiseries{m}$ denoted $d^{\star -1}$. The case $m=1$ of the next theorem is proved in \cite{Gray-KEF_SIAM2017}. 

\begth {\rm (Multiplicative dynamic output feedback group)}
\label{thm:mult_dyn_feedb_grp} 
$\left(\allpiseries{m},\star, \doubleone\right)$ forms a group with the identity element $\doubleone$. 
\endth
\medskip

The following statement is a direct interpretation of  Theorem~\ref{thm:mixprocomp_right_act_monoid} and the defintion of the group product $\star$.

\begth\label{th:pre-grp} $\left(\allpiseries{m}, \shuffle , \lhook\right)$ is a right pre-group with $\left(\allpiseries{m}, \star\right)$ as its Grossman--Larson group.
\endth
\medskip
 
Theorem~\ref{thm:mult_dyn_feedb_grp}, Theorem~\ref{thm:shuff_dist_mixprc} and \rref{eqn:grp_prod} imply the following relations for all $c \in \allpiseries{m}$:
\begin{align}
\label{eqn:inverse_relations}
\begin{aligned}
	c^{\star -1} 
		&= c^{\shuffle -1} \lhook c^{\star -1}\\
	\left(c^{\star -1}\right)^{\shuffle -1} 
		&= c \lhook c^{\star -1}. 
\end{aligned}
\end{align} 

\begle
\label{lem:mult_feedb_sub} 
Let $c,d \in \allpiseries{m}$, then $\left(c\star d\right)\left(\mathbf{1}\right) = c\left(\mathbf{1}\right)d\left(\mathbf{1}\right)$.
\endle

\begpr
Since $\left(c\lhook d\right) \left(\mathbf{1}\right) = c\left(\mathbf{1}\right)$, if follows that 
$$
	\left(c\star d\right)\left(\mathbf{1}\right) 
		= \left(d \shuffle \left(c\lhook d\right)\right)\left(\mathbf{1}\right)
		= \left(c\lhook d\right)\left(\mathbf{1}\right) d\left(\mathbf{1}\right)
		= c\left(\mathbf{1}\right)d\left(\mathbf{1}\right).
$$
\endpr

\begth
\label{thm:subgroup_mult_dyn_feedb} 
Let $M^{m} := \lbrace \doubleone + c \: : c \in \allproperseries{m}\rbrace$, then $\left( M^m, \star, \doubleone\right)$ is a normal subgroup of the multiplicative dynamic feedback group. 
\endth 

\begpr
Follows  from Lemma~\ref{lem:mult_feedb_sub}.
\endpr

 \smallskip
 
The group $M^m$ is isomorphic to the character group of the Hopf algebra $\mathcal{H}$ which is used for computation of feedback detailed in Section~\ref{sec:computational_framework}.


\section{Chen--Fliess Series Under Multiplicative Dynamic Output Feedback}
\label{sec:mult_dyn_feedb}

\begin{figure}[tb]
	\vspace*{-0.4 in}
	\begin{center}
		\includegraphics[scale = 0.7]{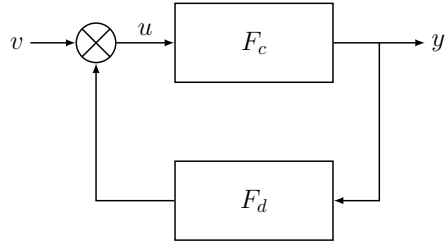}
		\caption{$F_c$ in multiplicative output feedback with $F_d$}
		\label{fig:mult_dyn_out_feedb}
	\end{center}
\end{figure}

Let $F_c$ be a Chen--Fliess series with $c \in \allseries{q}$. Suppose it is interconnected with $F_d$ where $d \in \allpiseriesXpri{m}$, as shown in Figure~\ref{fig:mult_dyn_out_feedb}. Note that, $\abs{X} = m+1$ and $\abs{X^{\prime}} = q+1$. 

Assume that the closed--loop system has a Chen--Fliess series representation, say $y = F_e[v]$, where $e \in \allseries{q}$. Then,
\begin{align*}
	y 
	= F_e[v] 
	= F_c[u] 
	= F_c[vF_d[y]] 
	= F_c[vF_d[F_e[v]]] 
	&= F_c[vF_{d\circ e}[v]]\\
	&= F_{c \lhook \left(d \circ e\right)}[v],
\end{align*}   
for any admissible input $v$. Hence, the series $e$ is characterized by the fixed point equation
\begin{align}
\label{eqn:dyn_feedb_fixed_point_eqn}
	e = c \lhook \left(d\circ e\right).
\end{align}

Theorem~\ref{thm:comp_strong_cont} and Theorem~\ref{thm:mixpr_cont} imply that $e \mapsto c \lhook \left(d \circ e \right)$ yields a strong contraction map in the ultrametric topology. Hence \rref{eqn:dyn_feedb_fixed_point_eqn} has a unique fixed point.
 
\begth
\label{thm:fixed_point_dyn_feedb} 
The series $c \lhook \left(d^{\shuffle -1}\circ c\right)^{\star -1} \in \allseries{q}$ is the unique fixed point of the map  $e \mapsto c \lhook \left(d \circ e \right)$ in the ultrametric space $\left(\allseries{q},\kappa\right)$.  
\endth

\begpr 
Let $e = c \lhook \left(d^{\shuffle -1}\circ c\right)^{\star -1}$, then using Theorem~\ref{thm:mix_assoc_comp_mixpr} and Theorem~\ref{thm:shuff_dist_mixprc}, gives
\begin{align*}
	c \lhook \left(d \circ e \right) 
	&= c \lhook \left[d \circ \left(c \lhook \left(d^{\shuffle -1}\circ c\right)^{\star -1}\right)\right]\\
	&= c \lhook \left[\left(d\circ c\right) \lhook \left(d^{\shuffle -1}\circ c\right)^{\star -1} \right]\\
	&= c \lhook \left[\left(d \circ c \right)^{\shuffle -1} \lhook \left(d^{\shuffle-1}\circ c\right)^{\star -1}\right]^{\shuffle -1}\\
	&\stackrel{\eqref{eq:distrib}}{=} 
		c \lhook \left[\left(d^{\shuffle -1}\circ c\right) \lhook \left(d^{\shuffle-1}\circ c\right)^{\star -1}\right]^{\shuffle -1}\\
	&\stackrel{\rref{eqn:inverse_relations} }{=} 
		c \lhook \left[\left(\left(d^{\shuffle -1}\circ c\right)^{\star -1}\right)^{\shuffle -1}\right]^{\shuffle -1} \\
	&= c \lhook \left(d^{\shuffle -1} \circ c\right)^{\star -1} 
	= e.
\end{align*}
\endpr

\begth
\label{thm:mult_dyn_feedb_formula} 
Given series $c \in \allseries{q}$ and $d \in \allpiseriesXpri{m}$ (such that $\abs{X} = m+1$ and $\abs{X^{\prime}} = q+1$), the generating series for the closed--loop system in Figure~\ref{fig:mult_dyn_out_feedb} is given by the {\em multiplicative dynamic feedback product} 
\begin{align}
\label{multidynfeedbprod} 
	c\check{@}d := c \lhook \left(d^{\shuffle -1}\circ c\right)^{\star -1}.
\end{align}
\endth

\medskip

The notion that feedback can be described in mathematical terms as a transformation group acting on the plant is well established in control theory \cite{Brockett_78}. The following theorem describes the situation in the present context.

\begth
\label{thm:mult_dyn_feedb_grp_action} 
The multiplicative dynamic feedback product \eqref{multidynfeedbprod} is a right group action by the group $\left(\allpiseriesXpri{m},\shuffle,\doubleone\right)$ on the set $\allseries{q}$, where $\abs{X} = m+1$ and $\abs{X^{\prime}} = q+1$.  
\endth

\begpr 
Let $c \in \allseries{q}$. Theorem~\ref{thm:mult_dyn_feedb_formula} gives
\begin{align*}
	c \check{@}\doubleone 
	&\stackrel{\eqref{multidynfeedbprod}}{=} c\lhook \left(\doubleone^{\shuffle -1} \circ c\right)^{\star-1}
	= c \lhook \doubleone = c.
\end{align*}
Let $d_1,d_2 \in \allpiseriesXpri{m}$, then using Theorem~\ref{thm:mixprocomp_right_act_monoid} together with the fact that the group inverse is an anti-homomorphism with respect to the group product, gives
\begin{align*}
	\left(c\check{@}d_1\right)\check{@}d_2 
	&= \left(c\check{@}d_1\right) \lhook \left(d_2^{\shuffle -1}\circ \left(c \check{@}d_1\right)\right)^{\star -1} \\
	&= \left(c \lhook \left(d_1^{\shuffle -1} \circ c\right)^{\star -1}\right) 
	\lhook \left(d_2^{\shuffle -1}\circ \left(c \lhook \left(d_1^{\shuffle -1}\circ c\right)^{\star -1}\right)\right)^{\star -1}\\
	&\stackrel{\eqref{eqn:rightact}}{=} \left(c \lhook \left(d_1^{\shuffle -1} \circ c\right)^{\star -1}\right) 
	\lhook \left(\left(d_2^{\shuffle -1}\circ c\right)\lhook \left(d_1^{\shuffle -1}\circ c\right)^{\star -1}\right)^{\star -1}\\
	&= c \lhook \Big[ \left(d_1^{\shuffle -1}\circ c\right)^{\star -1} \star \left(\left(d_2^{\shuffle -1}\circ c\right)
	\lhook \left(d_1^{\shuffle -1}\circ c\right)^{\star -1}\right)^{\star -1} \Big] \\
	&=c \lhook \Big[\left(\left(d_2^{\shuffle -1}\circ c\right) 
	\lhook \left(d_1^{\shuffle -1}\circ c\right)^{\star -1}\right) \star \left(d_1^{\shuffle -1}\circ c\right)\Big]^{\star -1}\\
	&\stackrel{\rref{eqn:grp_prod}}{=} c \lhook \Bigg[\left(d_1^{\shuffle -1} \circ c\right) \shuffle
	\Bigg(\Big(\left(d_{2}^{\shuffle -1}\circ c\right) \lhook \left(d_{1}^{\shuffle -1}\circ c\right)^{\star -1}\Big) 
	\lhook \left(d_1^{\shuffle -1}\circ c\right)\Bigg)\Bigg]^{\star-1}\\
	&= c \lhook \Bigg[\left(d_1^{\shuffle -1} \circ c\right) \shuffle
	\Bigg(\left(d_{2}^{\shuffle -1}\circ c\right) 
	\lhook \left(\left(d_1^{\shuffle -1}\circ c\right)^{\star -1} \star \left(d_1^{\shuffle -1}\circ c\right)\right)\Bigg)\Bigg]^{\star -1} \\
	&= c \lhook \left(\left(d_{1}^{\shuffle -1}\circ c\right) \shuffle \left(\left(d_2^{\shuffle -1} \circ c\right) 
	\lhook \doubleone\right)\right)^{\star -1} \\
	&= c \lhook \left(\left(d_{1}^{\shuffle -1}\circ c\right) \shuffle \left(d_{2}^{\shuffle -1}\circ c\right)\right)^{\star -1}\\
	&\stackrel{\eqref{eq:distrib}}{=} c \lhook \left(\left(d_1^{\shuffle -1}\shuffle d_2^{\shuffle -1}\right)\circ c\right)^{\star -1} \\
	&= c \lhook \left(\left(d_1\shuffle d_2\right)^{\shuffle -1}\circ c\right)^{\star -1}.
\end{align*}
Therefore, the statement is confirmed as the following identity holds
\begin{align*}
	\left(c\check{@}d_1\right)\check{@}d_2 
		&= c \check{@}\left(d_1 \shuffle d_2\right).
\end{align*}
\endpr
\medskip

\begre~The transformation group for the {\em additive dynamic feedback product} is the additive group $(\allseriesXpri{m},+,0)$ while for multiplicative feedback, it is $(\allpiseriesXpri{m},\shuffle,\doubleone)$.


\section{Invariance of Class and Relative Degree}
\label{sec:rel_degree_mult_feedb}

The notion of relative degree of a plant is essential in studying feedback linearization and system inversion~\cite{Isidori_95} as well as flatness~\cite{Levine_09}. The definitions of class and relative degree of a SISO Chen--Fliess series are characterized by the {\em structure} of its generating series \cite{Gray-etal_AUTO14,Gray_GS_relative degree,GS_thesis}. The definition of relative degree of a Chen--Fliess series is consistent with the classical definition whenever $y=F_c[u]$ has an input-affine analytic state space realization \cite{Gray-etal_AUTO14,Gray-KEF_SIAM2017}. 
\medskip

In this section, the alphabet $X = \{ x_0,x_1\}$ and $\supp(c)$ denotes the support of $c \in \allseries{}$. Recall that $c_N,c_F$ are the natural respectively forced part of the series $c$. Recall that we write $1$ for $1\emptyset$ in $ \allseries{}$.


\subsection{Class}
\label{ssec:class}

\begde
\label{de:class}~\cite{Gray_GS_relative degree}
A series $c\in\allseries{}$ is said to be of $r$--class, denoted by $\mathscr{C}(c) = r$, if $\supp(c_F) \subseteq x_{0}^{r-1}X^{+}$ and $\supp(c_F) \nsubseteq x_{0}^{r}X^{+}$. By definition, $\mathscr{C}(c) = \infty$ if $c_F = 0$.
\endde

\medskip

\begle~\cite{Gray_GS_relative degree} 
Every series $c\in \allseries{}$ has a class.
\endle

\medskip

\begex~\label{ex:class}
Let $c = 1 + x_0x_1^{2} + x_{0}^{2}x_{1}$. Thus, $c_F = x_0x_1^{2} + x_{0}^{2}x_{1}$, which implies $\supp(c_F) \subseteq x_{0}X^{+}$ but $\supp(c_F) \nsubseteq x_{0}^{2}X^{+}$. Thus, $\mathscr{C}(c) = 2$.
\endex
\medskip

The following lemma is essential in quantifying the class of $c \lhook d$ for all $c,d \in \allseries{}$. 

\begle
\label{lem:class_mixprocomp} 
Let $c,c^{\prime},d \in \allseries{}$ such that $\supp\left(c^{\prime}\right) \not\subseteq x_0\allwords$. Then 
$$
	i)\; x_0^k \lhook d = x_0^k \; \;\forall k \in \nat_0 
	\qquad	
	ii)\; c_N \lhook d = c_N
	\qquad
	iii)\; \supp\left(c^{\prime}\lhook d\right) \not\subseteq x_0\allwords.
$$
\endle

\begpr Statements $i)$ and $ii)$ are clear from \eqref{eqn:mixprocomp_def}.
\smallskip

Statement $iii)$ follows from $\supp\left(c^{\prime}\right) \not \subseteq x_0\allwords$ implying the existence of $x_i\eta \in \supp\left(c^{\prime}\right)$ where $x_i \neq x_0$ and $\eta\in \allwords$. Then, $x_i \eta \lhook d = x_i\left(d_i \shuffle \left(\eta\lhook d\right)\right)$. Therefore, $\supp\left(x_i\eta \lhook d\right) \subseteq x_iX^{\ast}$ implies that $ \supp\left(c^{\prime}\lhook d\right) \not \subseteq x_0X^{\ast}$.   
\endpr

\begth
\label{thm:class_mixprocomp} 
Let $c,d \in \allseries{}$, then the $r$-class $\mathscr{C}\left(c\lhook d\right) = \mathscr{C}\left(c\right)$.
\endth
 
 \begpr 
 Let $c \in \allseries{}$ be of $r$--class $\mathscr{C}\left(c\right)$, then  $c$ can be decomposed as $c = c_N +  x_0^{r-1}c^{\prime}$,
 where $c^{\prime}$ is a proper series such that $\supp\left(c^{\prime}\right) \not \subseteq x_0X^{\ast}$. Then,
 \begin{align*}
 	c \lhook d 
 	&= \left(c_N \lhook d\right) + x_0^{r-1}\left(c^{\prime}\lhook d\right) = c_N + x_0^{r-1}\left(c^{\prime}\lhook d\right), 
\end{align*}
with $\supp\left(c^{\prime}\lhook d\right) \not\subseteq x_0\allwords$. Hence, $\supp\left(c\lhook d\right)_F \subseteq x_0^{r-1}X^{+}$ and $\supp\left(c\lhook d \right)_F \not \subseteq x_0^{r}X^{+}$. 
\endpr

\begex Consider $c$ in Example~\ref{ex:class} and $d = 1+x_1 \in \allseries{}$. Then  $\mathscr{C}\left(c\lhook d\right) = 2 = \mathscr{C}\left(c\right)$ as $c \lhook d = 1 +x_0x_1^2 + 3x_0x_1^3 + 3x_0x_1^4 + x_0^2x_1 +x_0^2x_1^2.$
\endex 

The multiplicative dynamic feedback product preserves the class of its left argument.

\begth
\label{thm:class_invariance} 
Let $c \in \allseries{}$ be of $r$--class and $d \in \allpiseries{}$, then $\mathscr{C}\left(c\check{@}d\right) = r = \mathscr{C}\left(c\right)$.
\endth

\begpr Observe that $ \mathscr{C}(c\check{@}d) = \mathscr{C}( c \lhook \left(d^{\shuffle -1}\circ c\right)^{\star-1}) = r = \mathscr{C}\left(c\right). $
\endpr

\begex
\label{ex:class_invariance_dyn_feedb} 
Consider $c = x_1$ and $d = 1 + \sum_{k \in \nat} k!x_1^k$ in $\allseries{}$. Theorem~\ref{thm:mult_dyn_feedb_formula} implies $\mathscr{C}\left(c\check{@}d\right) = \mathscr{C}\left(c\right)=1$ as $c\check{@}d = x_1 + x_1x_0x_1 + 3x_1x_0x_1x_0x_1 + 4x_1x_0^2x_1^2 + \cdots.$
\endex


\subsection{Relative Degree}
\label{ssec:reldegree}

\begde \cite{Gray_GS_relative degree}
\label{de:relative-degree-c}
A series $c \in \allseries{}$ has relative degree $r$ if $\mathscr{C}(c) = r$ and the word $x_{0}^{r-1}x_{1} \in \supp(c_F)$. Otherwise, $c$ does not have relative degree.
\endde

The following theorem quantifies the relative degree of $c \lhook d$ for all $c,d \in \allseries{}$.

\begth
\label{thm:mixprocomp_rel_degree}
If $c \in \allseries{}$ with relative degree $r_c$ and $d \in \allseries{}$ such that $d\left(\mathbf{1}\right) \neq 0$, then $c \lhook d$ has relative degree $r_c$.
\endth

\begpr 
From Theorem~\ref{thm:class_mixprocomp}, $\mathscr{C}\left(c\lhook d\right) = r_c$. It remains to prove that $x_0^{r_c-1}x_1 \in \supp\left(c\lhook d\right)$. Since $c$ has relative degree $r_c$, it can be decomposed as: $c = c_N + \lambda x_0^{r_c-1}x_1 + x_0^{r_c-1}c^{\prime}$, where $\lambda \neq 0$ and $c^{\prime}$ is a proper series such that $x_1 \not \in \supp\left(c^{\prime}\right)$. Then,
\begin{align*}
	c \lhook d 
 	&= \left(c_N \lhook d\right) + \lambda \left(x_0^{r_c-1}x_1 \lhook d\right) + \left(x_0^{r_c-1}c^{\prime}\lhook d\right).
\end{align*}
Using~\rref{eqn:mixprocomp_def} and Lemma~\ref{lem:class_mixprocomp}, $c \lhook d = c_N + \lambda x_0^{r_c-1}x_1d + x_0^{r_c-1}\left(c^{\prime} \lhook d\right)$. Since $d = \alpha + d^{\prime}$, where $\alpha = d\left(\mathbf{1}\right) \neq 0$ and $d^{\prime}\left(\mathbf{1}\right) = 0$. Hence, 
\begin{align*}
	c \lhook d 
	= c_N + \lambda\alpha x_0^{r_c-1}x_1 + x_0^{r_c-1}x_1d^{\prime} 
	+ x_0^{r_c-1}\left(c^{\prime}\lhook d\right) \quad \left(\alpha\lambda \neq 0\right).
\end{align*}
Note that $x_1 \not \in \supp\left(c^{\prime}\right)$ implies $x_1 \not \in \supp\left(c^{\prime}\lhook d\right)$. 
\endpr
\medskip

\begex 
Let $c = 1 + x_0^2 + x_0x_1 + x_0^2x_1$ and $d = 1+x_1$. The relative degree $r_c = 2$ and $d\left(\mathbf{1}\right) = 1$. Then $r_{c \lhook d}= r_c$ as $c \lhook d = 1 + x_0^2 + x_0x_1 + x_0x_1^2 + x_0^2x_1 + x_0^2x_1^2.$
\endex

\begth
\label{thm:rel_degree_mult_feedb} 
If $c \in \allseries{}$ with relative degree $r_c$ and $d \in \allpiseries{}$, then $r_{c\check{@}d}=r_c$.
\endth

\begpr From Theorem~\ref{thm:mult_dyn_feedb_formula}, the multiplicative dynamic feedback product
\begin{align*}
	c \check{@} d = c\lhook \left(d^{\shuffle-1}\circ c\right)^{\star -1}.
\end{align*}
From $d \in \allpiseries{}$ being equivalent to $\left(d^{\shuffle-1}\circ c\right)^{\star -1} \in \allpiseries{}$ via Theorem~\ref{thm:proper_preser_comp} and Theorem~\ref{thm:mult_dyn_feedb_grp}, we deduce that $c \check{@} d$ has relative degree $r_c$ by Theorem~\ref{thm:mixprocomp_rel_degree}. 
\endpr
\medskip

\begex 
Consider Example~\ref{ex:class_invariance_dyn_feedb} where $r_c = 1$. Then $r_{c\check{@}d} = 1 = r_c$.
\endex


\section{Computational Framework}
\label{sec:computational_framework}

This subsection presents the Hopf algebra $\mathcal{H}$ of coordinate functions with character group isomorphic to the normal subgroup $M^m$ of multiplicative dynamic feedback group.  
\medskip

The vector space $V$ of coordinate maps on $\allseries{m}$ is spanned by $a^j_\eta$, where $\eta \in \allwords$ and $j = 1,2,\ldots,m$. For a series $c \in \allseries{m}$
\begin{align*}
	a^j_\eta \left(c\right) 
	=  c_j\left(\eta\right), \quad \forall\, j = 1,2,\ldots,m.
\end{align*}
The vector space $V = \bigoplus_{n \geq 0} V_n$ is graded, where $V_n$ is spanned by $a^j_\eta$, $j = 1,2,\ldots,m$, $|\eta|=n$.
\medskip

For all $k = 0,1,2,\ldots,m$ define a linear endomorphism $\theta_k: V \longrightarrow V$ such that $\theta_k(a^j_\eta)= a^j_{x_{k}\eta},$ for all $j = 1,2,\ldots,m$. Let $\mathcal{B}$ denote the graded symmetric algebra with grading induced by $V$. Its product is denoted by $\HAprod$ and the unit is $\mathsf{1}$. $\mathcal{B}$ is a bialgebra with the cocommutative coproduct $\Delta_{\shuffle}$ defined for all $c,d \in \allseries{m}$ as 
\begin{align*}
	\Delta_{\shuffle}(a^j_\eta)\left(c \otimes d\right) 
	&= \left(c\shuffle d\right)_{j}\left(\eta\right) 
	= \left(c_j \shuffle d_j\right)\left(\eta\right).
\end{align*}
By extending the usual unshuffle coproduct on words, $\Delta_{\shuffle}(\eta)= \sum_{(\eta)} {\eta'} \otimes {\eta''}$ (employing Sweedler's notation), it is understood that for all $a^j_\eta \in V$, 
\begin{align*}
	\Delta_{\shuffle}\left(a^j_\eta\right) 
	= \sum_{(\eta)} a^j_{\eta'} \otimes a^j_{\eta''}.
\end{align*}
The counit $\nu$ is defined as 
\begin{align*}
	\nu\left(h\right) 
	&= \begin{cases}
		1;\quad \text{if}\; h = \mathsf{1}, a^1_{\mathbf{1}}, a^2_{\mathbf{1}}, \ldots, a^m_{\mathbf{1}}\, \\
		0; \quad \text{otherwise}.
	\end{cases}
\end{align*}

\begth\cite{Foissy_15} 
On $V$
\begin{align*}
	\Delta_{\shuffle} \circ \theta_k 
	= \left(\theta_k \otimes \id
		+ \id \otimes \theta_k \right) \circ \Delta_{\shuffle}, 
\end{align*}
for all $k = 0,1,2,\ldots,m$.
\endth
\medskip

\begre~~Observe that $\left(\mathcal{B},\HAprod,\mathsf{1},\Delta_{\shuffle},\nu\right)$ is not a connected graded bialgebra as the elements $a^j_{\mathbf{1}}$, $j = 1,2,\ldots,m$, are group-like but not invertible.
\medskip

Denote $\mathfrak{s_i} := \left(a^i_{\mathbf{1}} - \mathsf{1}\right)$ for $i = 1,2,\ldots,m$. The ideal $\left(\mathfrak{s_1},\mathfrak{s_2},\ldots,\mathfrak{s_m}\right)$, is a bi-ideal. Define $\mathcal{H} = \mathcal{B}/\left(\mathfrak{s_1},\mathfrak{s_2},\ldots,\mathfrak{s_m}\right)$ and thus, 
 
\begth 
$\left(\mathcal{H},\HAprod,\mathsf{1},\Delta_{\shuffle},\nu\right)$ is a graded connected bialgebra. The character group of $\left(\mathcal{H},\Delta_{\shuffle},\nu\right)$ is isomorphic to the shuffle group $\left(M^m, \shuffle\right) \cong \underbrace{\left(M, \shuffle\right) \times \left(M, \shuffle\right) \times \cdots \times \left(M, \shuffle\right)}_{m\;\text{times}}$.
\endth
\medskip

The reminder of the section is to detail another coalgebra structure compatible with the graded augmented algebra of $\mathcal{H}$.


\subsection{Multiplicative Feedback Bialgebra}\label{subsec:feedback-bi}

Define an unital algebra map $\rho : \mathcal{B} \longrightarrow \mathcal{B} \otimes \mathcal{B}$  such that
\begin{align*}
	\rho\left(a^j_\eta\right)\left(c\otimes d\right) 
	&= \left(c \lhook d\right)_j\left(\eta\right) = \left(c_j \lhook d\right)\left(\eta\right), 
\end{align*}
for all $c,d\in \allseries{m}$. The map $\rho$ is not coassociative (compare with Theorem~\ref{thm:mixprocomp_right_act_monoid}). 

\begprop\label{prop:red_coprod_comp} For all $i = 0,1,2,\ldots,m$ and $j = 1,2,\ldots,m$;
\begen[(i)]
\item $\rho \left(a_{\mathbf{1}}^i \right) = a_{\mathbf{1}}^i \otimes \mathsf{1}$.
\medskip

\item $\rho \circ \theta_0 = \left(\theta_0 \otimes \id_{\mathcal{B}}\right) \circ \rho$.
\medskip

\item $\rho \circ \theta_k (a^j_\eta) 
	= \left(\theta_k \otimes \HAprod\right) \circ \left(\rho \otimes \id_{\mathcal{B}}\right) 
				\circ\sum_{(\eta)} a^j_{\eta'} \otimes a^k_{\eta''}$,
\enden
for all $j,k = 1,2,\ldots,m$ and $\eta \in \allwords$.
\endprop

\begpr 
Let $c \in \allseries{m}$, then 
\begin{align*}
	c &= c\left(\mathbf{1}\right)\emptyset + \sum_{j=0}^m x_j\left(x_j^{-1}\left(c\right)\right).
\end{align*}
\medskip

Thus for all $c,d \in \allseries{m}$,
\begin{align}
\label{eqn:mixprocomp_det}
	c \lhook d 
	&= c\left(\mathbf{1}\right)\emptyset + x_0\left(x_0^{-1}\left(c\right)\lhook d\right) 
	+ \sum_{j=1}^m x_j\left(d_j \shuffle \left(x_j^{-1}\left(c\right)\lhook d\right)\right) \nonumber.
\end{align}

\begen[(i)]
\item Let $c,d \in \allseries{m}$. From \rref{eqn:mixprocomp_def},
\begin{align*}
	\rho\left(a_{\mathbf{1}}^i\right)\left(c\otimes d\right) 
	= a_{\mathbf{1}}^i \otimes \mathsf{1}\left(c\otimes d\right).
\end{align*}
\medskip

\item Let $a^j_\eta \in V$. Then,
\begin{align*}
	\left(\rho \circ \theta_0\right)(a^j_\eta)\left(c\otimes d\right) 
	&= \left(c \lhook d\right)_j\left(x_0\eta\right) 
	=\left(x_0^{-1}\left(c\lhook d\right)_j\right)\left(\eta\right). \\
\end{align*}
From \rref{eqn:mixprocomp_def},
\begin{align*}
 	\left(\rho \circ \theta_0\right)(a^j_\eta)\left(c\otimes d\right)
	&= \left(x_0^{-1}\left(c_j\right)\lhook d\right)\left(\eta\right) \\
	&= \rho\left(a^j_\eta\right)\left(x_0^{-1}\left(c\right)\otimes d\right) \\
	&= \left(\theta_0 \otimes \id_{\mathcal{H}} \right) \circ \rho\left(a^j_\eta\right)\left(c\otimes d\right).
\end{align*}
\medskip

\item For $k = 1,2, \ldots,m$;
\begin{align*} 
	\left(\rho \circ \theta_k\right)\left(a^j_\eta\right)\left(c\otimes d\right) 
	&= \left(\left(c\lhook d\right)_j\right)\left(x_k\eta\right)
	=\left(x_k^{-1}\left(c\lhook d\right)_j\right)\left(\eta\right).\\
\end{align*} 
From \rref{eqn:mixprocomp_def},
\begin{align*}
	\left(\rho \circ \theta_k\right)\left(a^j_\eta\right)\left(c\otimes d\right) 
	&= \left(d_i \shuffle x_k^{-1}\left(c_j\right)\lhook d\right) \left(\eta\right)\\
	&= \Big( \sum_{(\eta)} a^j_{\eta'} \otimes a^k_{\eta''} \Big)
	\left(\left(x_k^{-1} \left(c\right)\lhook d\right)\otimes d\right) \\
	&= \left(\rho\otimes \id_{\mathcal{B}}\right)\circ
	\Big( \sum_{(\eta)} a^j_{\eta'} \otimes a^k_{\eta''} \Big) \left(x_k^{-1}\left(c\right)\otimes d \otimes d\right) \\
	&= \left(\theta_k \otimes \HAprod\right)\circ \left(\rho \otimes \id_{\mathcal B}\right) 
	 \circ \Big( \sum_{(\eta)} a^j_{\eta'} \otimes a^k_{\eta''} \Big)\left(c\otimes d\right).
\end{align*}
\enden
\endpr

\medskip

Define the coproduct $\Delta: \mathcal{B} \longrightarrow \mathcal{B} \otimes \mathcal{B}$ such that 
\begin{align*}
\Delta\left(a^j_\eta\right)\left(c\otimes d\right) 
&= \left(c\star d\right)_{j}\left(\eta\right),
\end{align*}
for all $c,d \in \allseries{m}$. Since the product, $\star$, is associative~(see Theorem~\ref{thm:mult_comp_assoc}), the map $\Delta$ is coassociative.

\begprop\label{prop:red_coprod_hopf_sec} 
The coproduct on $V$ is defined as
\begin{align*}
\Delta = \left(\id_{\mathcal{B}} \otimes \HAprod\right) \circ 
\left(\rho \otimes \id_{\mathcal{B}} \right) \circ \Delta_{\shuffle}.
\end{align*}
\endprop 

\begpr 
Let $c,d \in \allseries{m}$ and $\eta \in \allwords$,  
\begin{align*}
\Delta\left(a^k_\eta\right) \left(c\otimes d\right) 
&= \left(c\star d\right)_k \left(\eta\right).
\end{align*}
Using \rref{eqn:grp_prod},
\begin{align*}
\Delta\left(a^k_\eta\right) \left(c\otimes d\right) 
&= \left(d_k \shuffle \left(c_k\lhook d\right)\right) \left(\eta\right) \\
&= \Delta_{\shuffle}\left(a^k_\eta\right)\left(\left(c \lhook d\right)\otimes d\right)\\
&=\left(\rho \otimes \id_{\mathcal{B}} \right)\circ \Delta_{\shuffle}\left(a^k_\eta\right)\left(c\otimes d\otimes d\right) \\
&= \left(\id_{\mathcal{B}} \otimes \HAprod\right) \circ 
\left(\rho \otimes \id_{\mathcal{B}} \right) \circ \Delta_{\shuffle}\left(a^k_\eta\right)\left(c\otimes d\right).
\end{align*}
\endpr
\medskip

\begre{Observe that $\left(\mathcal{B}, \HAprod, \mathsf{1},\Delta,\nu\right)$ is not a Hopf algebra as the elements $a_\mathbf{1}^j, j = 1,2\ldots,m$ are group--like but not invertible.}
\smallskip

The coproduct $\Delta$ respects the grading of $\mathcal{B}$. The following lemma is necessary to prove this claim. Recall that the length of a word $x_{i_1}x_{i_2}\cdots x_{i_k}$ is denoted $\abs{x_{i_1}x_{i_2}\cdots x_{i_k}} (= k)$. 

\begle
\label{lem:grade_preservation_B} 
For all $n \geq 0$;
\begin{align*}
\rho\left(V_n\right) \subseteq \bigoplus\limits_{i+j = n} V_i \otimes \mathcal{B}_{j} 
\end{align*}
\endle

\begpr 
The linear endomorphism $\theta_i$ is a homogeneous operator of degree $1$ on $V$ for all $i = 0,1,2,\ldots,m$.
\begin{align*}
|\theta_i a_{\eta}^j|
= \abs{x_i\eta} = \abs{\eta} + 1.
\end{align*} 
Therefore, $\theta_j : V_n \longrightarrow V_{n+1} $ for all $j = 1,\ldots,m$ and $n \geq 0$.
\medskip

The lemma is proven by induction on degree $n$. The base case $(n = 0)$ is true by Proposition~\ref{prop:red_coprod_comp}.

Assume that the statement holds for all $n \leq k$. Let $a_{\tilde\eta}^j \in V_{k+1}$. Then there are two possibilities:
\begen[(i)]
\item $a_{\tilde\eta}^j = a_{x_0 \eta}^j$ where $a_{\eta}^j \in V_k$. 
\begin{align*}
\rho \left(a_{\tilde\eta}^j \right) 
&= \left(\rho \circ \theta_0\right) \left(a_{\eta}^j\right) \\
&= \left(\theta_0 \otimes \id_{\mathcal{B}}\right) \circ \rho\left(a_{\eta}^j\right). 
\quad (\text{c.f.~Propostion}~\ref{prop:red_coprod_comp})
\end{align*}
Note that $\rho\left(a_{\eta}^j\right) \subseteq \bigoplus\limits_{i+j = k} V_i \otimes \mathcal{B}_j$ via induction hypothesis. Hence,
\begin{align*}
\rho\left(a_{\tilde\eta}^j \right) 
&\subseteq \bigoplus\limits_{i+j =k} \theta_0\left(V_i\right) \otimes \mathcal{B}_j
\subseteq \bigoplus\limits_{i+j = k+1} V_i \otimes \mathcal{B}_j.
\end{align*}
\medskip

\item $a_{\tilde\eta}^j = a_{x_p\eta}^j$ where $p \neq 0$ and $ a_{\eta}^j \in V_k$.
\begin{align*}
\rho\left(a_{\tilde\eta}^j \right) 
&= \left(\rho \circ \theta_p\right)\left(a_{\eta}^j\right) \\
&= \left(\theta_p \otimes \HAprod\right) \circ \left(\rho \otimes \id_{\mathcal{B}}\right) 
\circ\Big( \sum_{(\eta)} a^j_{\eta'} \otimes a^p_{\eta''} \Big) 
\quad \quad  (\text{c.f.~Proposition}~\ref{prop:red_coprod_comp}).
\end{align*}
Recall that $\sum_{(\eta)} a^j_{\eta'} \otimes a^p_{\eta''} \in \left(V \otimes V\right)_k$
\begin{align*}
\rho\left(a_{\tilde\eta}^j \right) 
&\subseteq \left(\theta_p \otimes \HAprod\right) \left(\bigoplus\limits_{i+j=k}\rho\left(V_i\right) 
\otimes \mathcal{B}_j\right)\\
&\subseteq \bigoplus\limits_{i_i+i_2+j = k} \theta_p\left(V_{i_1}\right) \otimes \HAprod\left(\mathcal{H}_{i_2} 
\otimes \mathcal{B}_j\right) 
\subseteq \bigoplus\limits_{i+j = k+1}V_{i} \otimes \mathcal{B}_j.
\end{align*}
\enden
\endpr

\begprop
\label{prop:grading-well-defined-B} 
For all $n \geq 0$;
\begin{align*}
\Delta\left(\mathcal{B}_n\right) \subseteq \bigoplus\limits_{i + j =n} \mathcal{B}_i \otimes \mathcal{B}_j
\end{align*}
\endprop

\begpr 
Since $\Delta$ is an algebra morphism, it is enough to prove that for all $n \geq 0$
\begin{align*}
\Delta\left(V_n\right) \subseteq \bigoplus\limits_{i + j = n} V_i \otimes \mathcal{B}_j.
\end{align*}
Hence, 
\begin{align*}
\Delta\left(V_n\right) 
&= \left(\id_{\mathcal B} \otimes \HAprod\right) \circ \left(\rho \otimes \id_{\mathcal B}\right) 
\circ \Delta_{\shuffle}\left(V_n\right) \quad \left(\text{c.f.~Proposition}~\ref{prop:red_coprod_hopf_sec}\right) \\
&\subseteq \left(\id_{\mathcal B} \otimes \HAprod\right) \circ \left(\rho \otimes 
\id_{\mathcal B}\right) \left(\bigoplus\limits_{i+j =n} V_i \otimes V_j\right) \\
&\subseteq \left(\id_{\mathcal B} \otimes \HAprod\right)\left(\bigoplus\limits_{i_1+i_2+j = n}
V_{i_1} \otimes \mathcal{B}_{i_2} \otimes V_j\right) \quad \left(\text{c.f.~Lemma}~\ref{lem:grade_preservation}\right) \\
&\subseteq \bigoplus\limits_{i+j = n} V_i \otimes \mathcal{B}_j
\end{align*}
\endpr
\medskip

The cointeraction between the bialgebras $\left(\mathcal{B},\Delta\right)$ and $\left(\mathcal{B}, \Delta_{\shuffle}\right)$ is established next.
 
\begth\label{th:bi-coaction} $\left(\mathcal{B},\HAprod,\mathsf{1},\Delta_{\shuffle},\nu\right)$ is a right graded comodule bialgebra of $\left(\mathcal{B},\HAprod,\mathsf{1},\Delta,\nu\right)$ with the coaction map $\rho$.
\endth

\begpr \begen[(i)]
\item $\rho$ is a coaction follows from its dual statements Theorem~\ref{thm:shuff_dist_mixprc} and Theorem~\ref{thm:mixprocomp_right_act_monoid}. 
\medskip

\item $\rho$ respects grading of $\mathcal{B}$ from Lemma~\ref{lem:grade_preservation_B}.
\enden 

\medskip

\begre~ $\left(\mathcal{H},\HAprod,\mathsf{1},\Delta_{\shuffle},\Delta,\nu,\rho\right)$ is a graded connected cointeracting bialgebras with the right coaction map $\rho$. For more details on cointeracting bialgebas and double bialgebras we refer to~\cite{Foissy_22}.
\medskip

\subsection{From cointeracting bialgebra $\mathcal{B}$ to cointeracting Hopf algebra $\mathcal{H}$}\hfill\smallskip

Observe that $\rho(a_{\mathbf{1}}^i - \mathsf{1}) = \left(a_{\mathbf{1}}^i - \mathsf{1}\right) \otimes \mathsf{1}$ for all $i = 1,2,\ldots,m$ and the coproduct 
\begin{align*}
\Delta\left(a_{\mathbf{1}}^{i} - \mathsf{1}\right) &= \left(a_{\mathbf{1}}^i \otimes a_{\mathbf{1}}^i\right) - \left(\mathsf{1} \otimes \mathsf{1}\right) \\
&=a_{\mathbf{1}}^i \otimes \left(a_{\mathbf{1}}^i - \mathsf{1}\right) + \left(a_{\mathbf{1}}^i - \mathsf{1}\right) \otimes \mathsf{1}. 
\end{align*}

Therefore, (as defined earlier in Section~\ref{sec:computational_framework}) for $\mathfrak{s_i} = a_{\mathbf{1}}^i-1$ and the ideal $\left(\mathfrak{s_1}, \mathfrak{s_2}, \ldots, \mathfrak{s_m}\right)$ is a bi--ideal of the bialgebra $\left(\mathcal{B}, \Delta\right)$. 
\smallskip

Defining $\mathcal{H} = \mathcal{B}/\left(\mathfrak{s_1},\mathfrak{s_2}, \ldots, \mathfrak{s_m}\right)$ provides a connected structure with $\mathcal{H}_0 \cong \re\mathsf{1}$, thus making $\left(\mathcal{H},\Delta\right)$ a Hopf algebra.

\medskip

\begre{ Let $\mathcal{I} := \left(\mathfrak{s_1},\mathfrak{s_2}, \ldots, \mathfrak{s_m}\right)$ denote the bi--ideal. Observe that $I$ is the bi--ideal for both the cointeracting bialgebras $\left(\mathcal{B}, \Delta_{\shuffle}\right)$ and $\left(\mathcal{B}, \Delta\right)$. Since,	 
\begin{align*}
\rho(\mathfrak{s_i}) \subseteq \mathfrak{s_i} \otimes \mathcal{B}.
\end{align*}
Thus, 
\begin{align*}
\rho\left(\mathcal{B}/\mathcal{I}\right) = \rho(\mathcal{H}) &\subseteq \rho(\mathcal{B}/\mathcal{I}) \otimes \mathcal{B} \\
&= \mathcal{H} \otimes \mathcal{B} \\
&= \mathcal{H} \otimes \left(\mathcal{H} + \mathcal{I}\right) \\
&= \left(\mathcal{H} \otimes \mathcal{H}\right) + \cancel{\left(\mathcal{H} \otimes \mathcal{I}\right)}.
\end{align*}

Therefore, $\rho : \mathcal{H} \longrightarrow \mathcal{H} \otimes \mathcal{H}$ is a right coaction map on Hopf algebra $\left(\mathcal{H},\Delta_{\shuffle}\right)$ by the Hopf algebra $\left(\mathcal{H},\Delta\right)$. 
}
\medskip

The analogous statements of coproduct $\Delta$ and the map $\rho$ respects the grading of $\mathcal{H}$ follows from Lemma~\ref{lem:grade_preservation_B} and Proposition~\ref{prop:grading-well-defined-B} respectively.

\begle \label{lem:grade_preservation} 
For all $n \geq 1$;
\begin{align*}
	\rho\left(V_n\right) \subseteq \bigoplus\limits_{i+j = n} V_i \otimes \mathcal{H}_{j} 
\end{align*}
\endle

\medskip
\begth\label{prop:grading-well-defined-H} \hfill  
\begen
\item For all $n \geq 0$;
\begin{align*}
\Delta\left(\mathcal{H}_n\right) \subseteq \bigoplus\limits_{i + j =n} \mathcal{H}_i \otimes \mathcal{H}_j. 
\end{align*}
\medskip

\item The graded Hopf algebra $\mathcal{H}$ is right--handed viz. for all $n \geq 1$
\begin{align*}
\Delta'\left(V_n\right) \subseteq \bigoplus_{\substack{i+j=n \\i,j \geq 1}} V_{i} \otimes \mathcal{H}_{j} 
\end{align*}
where $\Delta'(\cdot) = \Delta(\cdot) - \mathsf{1} \otimes \cdot - \cdot \otimes \mathsf{1}$ is the reduced coproduct.
\enden
\endth
\medskip

\begth\label{thm:graded_Hopf} 
$\left(\mathcal{H},\HAprod,\mathsf{1},\Delta,\nu\right)$ is a graded connected (right--handed) bialgebra. The character group of $\left(\mathcal{H},\Delta,\nu\right)$ is isomorphic to the multiplicative feedback group $\left(M^m, \star\right)$ and
\begin{align*}
\left(M^m, \star\right) \not\cong \underbrace{\left(M,\star\right) \times \left(M,\star\right) \times 
	\cdots \times \left(M, \star\right)}_{m~\text{times}} \quad \quad \text{for}\; m > 1.
\end{align*}
\endth
\medskip

The summary of this section so far:
\begen
\item $\left(\mathcal{H},\HAprod,\mathsf{1},\Delta_{\shuffle},\nu\right)$ is a graded connected bialgebra. 
\medskip

\item The character group~of $\left(\mathcal{H},\Delta_{\shuffle},\nu\right)$ is isomorphic to the group $\left(M^m,\shuffle\right)$.
\medskip

\item $\left(\mathcal{H},\HAprod,\mathsf{1},\Delta,\nu\right)$ is a graded connected right--handed bialgebra. 
\medskip

\item The character group~of $\left(\mathcal{H},\Delta,\nu\right)$ is isomorphic to the group $\left(M^m,\star\right)$.
\enden
\medskip

The following results stem from Theorem~\ref{th:bi-coaction}.
\begth\label{th:comodule-str} $\left(\mathcal{H},\HAprod,\mathsf{1},\Delta_{\shuffle},\nu\right)$ is a right graded comodule Hopf algebra of $\left(\mathcal{H},\HAprod,\mathsf{1},\Delta,\nu\right)$ with the (graded) coaction map $\rho$.
\endth
\medskip

\begre{~Theorem~\ref{th:comodule-str}, along with the definition of coproduct $\Delta$ as $\Delta = \left(\id \otimes \HAprod\right) \circ \left(\rho \otimes \id\right) \circ \Delta_{\shuffle}$ implies that the character group $\left(M^m, \star\right)$ is the Grossman--Larson group of the pre-group $\left(M^m, \shuffle, \lhook\right)$ as seen in Theorem~\ref{th:pre-grp}}. 
\medskip

A few examples of the computation of $\rho$ on $\mathcal{H}$ are (for all $i,j,k = 1,2,\ldots,m$):

\begin{align*}
\rho \left(a_{x_0}^i\right) 
&= a_{x_0}^i \otimes \mathsf{1} 		\\
\rho\left(a_{x_i}^j\right) 
&=  a_{x_i}^j\otimes \mathsf{1}			 \\
\rho \left(a_{x_0x_0}^i\right) 
&=a_{x_0x_0}^i \otimes \mathsf{1}		 \\
\rho\left(a_{x_0x_i}^j\right) 
&= a_{x_0x_i}^j \otimes \mathsf{1}		 \\
\rho\left(a_{x_ix_0}^j \right) 
&= \left(a_{x_ix_0}^j \otimes \mathsf{1}\right) 
+ \left(a_{x_i}^j \otimes a_{x_0}^i\right)	 \\
\rho\left(a_{x_ix_j}^k\right) 
&= \left(a_{x_ix_j}^k\otimes \mathsf{1}\right) + \left(a_{x_i}^k \otimes a_{x_j}^i\right). 
\end{align*}
\medskip

A few examples of the coproduct $\Delta$ computation are:
\begin{align*}
\Delta\left(a_{x_0}^i\right) 
&= \left(a_{x_0}^i\otimes \mathsf{1}\right) 
+ \left(\mathsf{1} \otimes a_{x_0}^i\right)		\\
\Delta\left(a_{x_i}^j\right)  
&= \left(a_{x_i}^j \otimes \mathsf{1}\right) 
+ \left(\mathsf{1} \otimes a_{x_i}^j\right)		\\
\Delta\left(a_{x_0x_0}^i\right) 
&= \left(a_{x_0x_0}^i \otimes \mathsf{1}\right) 
+ 2\left(a_{x_0}^i \otimes a_{x_0}^i\right) 
+ \left(\mathsf{1} \otimes a_{x_0x_0}^i\right)	\\
\Delta\left(a_{x_0x_i}^j\right) 
&= \left(a_{x_0x_i}^j \otimes \mathsf{1}\right) 
+ \left(a_{x_0}^j \otimes a_{x_i}^j\right) 
+ \left(a_{x_i}^j \otimes a_{x_0}^j\right) 
+ \left(\mathsf{1} \otimes a_{x_0x_i}^j \right)	\\
\Delta\left(a_{x_ix_0}^j\right) 
&= \left(a_{x_ix_0}^j \otimes \mathsf{1}\right) 
+ \left(a_{x_i}^j \otimes a_{x_0}^i \right) 
+ \left(a_{x_i}^j \otimes a_{x_0}^j \right) 
+ \left(a_{x_0}^j \otimes a_{x_i}^j \right) 
+ \left(\mathsf{1} \otimes a_{x_ix_0}^j \right)	\\
\Delta\left(a_{x_ix_j}^k\right) 
&= \left(a_{x_ix_j}^k \otimes \mathsf{1}\right) 
+ \left(a_{x_i}^k \otimes a_{x_j}^i \right) 
+ \left(a_{x_i}^k\otimes a_{x_j}^k \right) 
+ \left(a_{x_j}^k \otimes a_{x_i}^k\right) 
+ \left(\mathsf{1} \otimes a_{x_ix_j}^k\right).
\end{align*}

\begco
\label{cor:coprod_for_x_0} 
If $n \in \nat_0$, then for all $i = 1,2,\ldots,m$;
\begin{align*}
\rho\left(a_{x_0^n}^i \right) 
&= a_{x_0^n}^i \otimes \mathsf{1}. \\ 
\Delta\left(a_{x_0^n}^i\right) 
&= \sum\limits_{k = 0}^n \binom{n}{k} a_{x_0^k}^i \otimes a_{x_0^{n-k}}^i,
\end{align*}
where $a_{x_0^0}^i:= \mathsf{1}$.
\endco

\begpr 
The proof of is by induction on $n$. 
\smallskip
The second identity follows from Proposition~\ref{prop:red_coprod_hopf_sec} and \rref{eqn:shuffle_pow}.
\endpr
\medskip 

The Example~$4.10$ in~\cite{Gray-KEF_SIAM2017} for the computation of feedback group inverse of a one-dimensional series is reworked within our framework.
\medskip

\begex 
Let $c = 1 - x_1 \, \in \re[[x_1]]$. Let $S_{\star}$ denote the antipode of the graded connected Hopf algebra $\left(\mathcal{H},\HAprod,\mathsf{1},\Delta,\nu\right)$.
 
\begin{align*}
	\left(a_{x_1}^1\right)(c^{\star -1}) 
	= S_\star\left(a_{x_1}^1\right)\left(c\right) 
	= -\left(a_{x_1}^1\right)(c) 
	= 1.
\end{align*} 
The reduced coproduct $\Delta^{\prime}\left(\cdot\right) := \Delta\left(\cdot\right) - \mathsf{1} \otimes \cdot - \cdot \otimes \mathsf{1}$. Hence
\begin{align*}
	\Delta^{\prime}\left(a_{x_1x_1}^1\right) 
	= 3a_{x_1}^1 \otimes a_{x_1}^1.
\end{align*}
Therefore,
\begin{align*}
	\left(a_{x_1}^1\right)\left(c^{\star -1}\right) 
	&= S_\star \left(a_{x_1x_1}^1\right)\left(c\right) 
		= -\left(a_{x_1x_1}^1\right)\left(c\right) -
	 	\left(3\left(a_{x_1}^1\right).
		S\left(a_{x_1}^1\right)\right)\left(c\right)\\
	& = \left(-a_{x_1x_1}^1 + 3\left(a_{x_1}^1\right)^2\right)\left(c\right) = 3.
\end{align*}
Similarly for computing $a_{x_1x_1x_1}^1\left(c^{\star -1}\right)$;
\begin{align*}
	\Delta^{\prime}\left(a_{x_1x_1x_1}^1\right) 
	= \left(4a_{x_1}^1 \otimes a_{x_1x_1}^1\right) 
	+ \left(6a_{x_1x_1}^1 \otimes a_{x_1}^1\right) 
	+ 3 \left(a_{x_1}^1 \otimes \left(a_{x_1}^1\right)^2\right).
\end{align*}
Thus,
\begin{align*}
	\left(a_{x_1x_1x_1}^1\right)\left(c^{\star -1}\right) 
	&= \left[-a_{x_1x_1x_1}^1 - 4a_{x_1}^1 
	.S\left(a_{x_1x_1}^1\right) 
	-6a_{x_1x_1}^1
	.S\left(a_{x_1}^1\right) 
		- 3a_{x_1}^1.\left(S_\star \left(a_{x_1}^1\right)\right)^2\right](c)\\
	&= 0 -4(-1)(3) -6(0)(-1) -3 (-1)(1)^2 = 15.
\end{align*}

Continuing the computation one obtains $c^{\star -1} = 1 + x_1 + 3x_1^2 + 15x_1^3 + 105 x_1^4 + \cdots$. 
\endex

\begre{~For further part of the document, if $\left(\mathcal{C},\Delta, \epsilon\right)$ is a graded connected coalgebra then, $\Delta'$ denotes the reduced coproduct where 
\begin{align*}
\Delta' \left(\cdot\right) = \Delta\left(\cdot\right) - \cdot \otimes \mathsf{1} - \mathsf{1} \otimes \cdot
\end{align*}
}


\section{pre-Lie Structure on $\allproperpoly{m}$}\label{sec:pre-Lie}

Let $k$ be a field of characteristic zero and $V$ be a $k$--vector space. The definitions of (right) pre-Lie algebra and pre-Lie coalgebra are recollected briefly. For more details on the pre-Lie algebras, we refer to~\cite{CatPat2021,Manchon_11}. 

\begde\label{de:pre-Lie} $\left(V,\bullet\right)$ is a (right) pre-Lie algebra if the magmatic map $\bullet: V^{\otimes 2} \longrightarrow V$ satisfies for all $a,b,c \in V$ the (right) pre-Lie identity
\begin{align*}
	\left(a \bullet b\right) \bullet c - a \bullet \left(b \bullet c\right) 
		= \left(a \bullet c\right) \bullet b - a \bullet \left(c \bullet b\right).
\end{align*}
\endde

Define $\left[a,b\right]_{\bullet} := a\bullet b - b \bullet a$ , then $\left(V, \left[\cdot,\cdot\right]_{\bullet}\right)$ is a Lie algebra.
\medskip

\begde\label{de:copre-Lie} $(V,\mathring{\delta})$ is a (right) pre-Lie coalgebra if $\mathring{\delta}: V \longrightarrow V^{\otimes 2} $ satisfies
\begin{align*}
	(\id_V \otimes \id_V \otimes \id_V - \tau_{(23)}) \circ \Big(\big(\mathring{\delta} \otimes \id_V\big) - \big(\id_V \otimes \mathring{\delta}\big)\Big) \circ \mathring{\delta} 
	 = 0,
\end{align*}
where $\tau_{(23)}: V^{\otimes 3} \to V^{\otimes 3}$, $\tau_{(23)} (a \otimes b \otimes c) = a \otimes c \otimes b$.
\endde
\medskip

Let $\mathcal{S}(V)$ be the free symmetric algebra generated by the vector space $V$, with $\HAprod$ denoting the symmetric product.
\smallskip

\begprop
\label{prop: copre-Lie} 
Let $(\mathcal{S}(V),\HAprod, \delta,\Delta,\epsilon,\rho)$ be graded connected cointeracting bialgebra where 

\begen[(i)]
\item $\left(\mathcal{S}(V),\HAprod,\delta,\epsilon\right)$ is a graded connected Hopf algebra in the category of $\left(\mathcal{S}(V), \HAprod, \Delta, \epsilon\right)$ right comodule with coaction map $\rho$.

\item $ \Delta = \left(\id \otimes \HAprod\right) \circ \left(\rho \otimes \id \right) \circ \delta$.

\item $\delta'(V) \subseteq \mathcal V \otimes \mathcal{S}^{+}(V)$. 

\item For all $x \in V$: $\rho'(x) := \rho(x) - x \otimes 1 \subseteq V \otimes \mathcal{S}^{+}(V)$. 
\enden
Then,
\begen
\item $\left(\mathcal{S}(V),\HAprod,\Delta,\epsilon\right)$ is a right--handed bialgebra.
\medskip
\item On $V$: $\mathring{\Delta} = \mathring{\rho} + \mathring{\delta}$, where
$\mathring{\upsilon} := \left(\pi_V \otimes \pi_V\right) \circ \upsilon$ for all $\upsilon \in Hom\left(\mathcal{S}(V),\mathcal{S}(V)\otimes \mathcal{S}(V)\right)$.
\enden

where $\pi_V$ is the canonical projection onto $V$.
\endprop

\begpr 
\begen 
\item It is straightforward to check that bialgebra $\left(\mathcal{S}(V),\HAprod,\Delta\right)$ is right--handed.
\medskip

\item  On the vector space $V$, 
\begin{align*}
\mathring{\Delta} 
&= \left(\pi_V \otimes \pi_V\right) \circ \left(\id \otimes \HAprod\right) \circ \left(\rho \otimes \id\right) \circ \delta \\
&= \left(\id \otimes \pi_V\right) \circ \left(\pi_V \otimes \HAprod\right) \circ \lbrace\left[ \rho' + \left(\id \otimes 1\right)\right] \otimes \id\rbrace \circ \delta\\
&= \lbrace\left(\id \otimes \pi_V\right) \circ \left(\pi_V \otimes V\right) \circ \left(\rho' \otimes \id\right) \circ \delta\rbrace + \left(\pi_V \otimes \pi_V \right) \circ \delta \\
&= \lbrace \left(\id \otimes \pi_V\right) \circ \left(\pi_V \otimes \HAprod\right) \circ \left(\rho' \otimes \id\right) \circ \lbrace1\otimes\id + \id \otimes 1 + \delta'\rbrace\rbrace + \mathring{\delta}\\
&= \lbrace\left(\pi_V \otimes \pi_V\right) \circ \rho'\rbrace + \lbrace\left(\pi_V \otimes \pi_1\right) \circ \left(\id \otimes \HAprod\right) \circ \left(\rho' \otimes \id\right) \circ \delta'\rbrace + \mathring {\delta} \\
&= \mathring{\rho} + \mathring{\delta},
\end{align*}
since $\left(\pi_V \otimes \pi_V \right) \circ \rho' = \mathring{\rho} = \left(\pi_V \otimes \pi_V \right)$ and $\left(\id \otimes \HAprod\right) \circ \left(\rho' \otimes \id \right) \circ \delta' \in \ker\left(\pi_V \otimes \pi_V\right)$.
\enden
\endpr

\begre{~$\mathcal{S}(V)$ is a right--hand bialgebra with both coproducts, $\Delta$ and $\delta$, thus $V$ inherits a graded right pre-Lie coalgebra with $\mathring{\Delta}$ and $\mathring{\delta}$, as defined in Proposition~\ref{prop: copre-Lie}~\cite{MenousPatras2015}.}
\medskip


\subsection{pre-Lie Coproduct on $V^{+}$}
\label{subsec:copre-Lie on V}

Observe that from Section~\ref{sec:computational_framework}, the Hopf algebra $\mathcal{H} \cong \mathcal{S}(V^{+})$ as $\re$-algebras where $V^{+} = \bigoplus_{n \geq 1} V_n$.

\begco $V^{+}$ is a graded right pre-Lie coalgebra with the pre-Lie coproducts $\mathring{\Delta}$ and $\mathring{\Delta}_{\shuffle}$ with
\begeq
\label{eq:copre-Lie-rel}
	\mathring{\Delta} = \mathring{\rho} + \mathring{\Delta}_{\shuffle}.
\endeq 
\endco
\smallskip

\begre{~ It is straightforward to check that for all $a_{\eta}^{i} \in V^{+}$ 
\begin{align*}
\mathring{\Delta}_{\shuffle}\left(a_{\eta}^{i}\right) 
	&= \left(\pi_{V^{+}} \otimes \pi_{V^{+}}\right) \circ \Delta_{\shuffle}\left(a_{\eta^{i}}\right) \\
	&= \sum\limits_{\left(\eta\right)} a_{\eta'}^i \otimes a_{\eta''}^i
\end{align*}
where $\Delta_{\shuffle}'(\eta) = \sum\limits_{\left(\eta\right)} \eta' \otimes \eta''$ is the reduced unshuffle coproduct in $\allpoly$. 
}

The following proposition computes $\mathring{\rho}$ based on the induction on degree of $a_{\eta}^j \in V^{+}$.

\begth\label{th:rho-lin} For all $k,j = 1,2,\ldots,m$ and $i = 0,1,\ldots,m$
\begen[(i)] 
\item $\mathring{\rho}\left(a_{x_{i}}^k\right) = 0$.

\item $\mathring{\rho}\circ \theta_0 = \left(\theta_0 \otimes \id\right) \circ \mathring{\rho}$.

\item $\mathring{\rho}\circ\theta_k \left(a_{\eta}^j\right) = \left(\theta_k \otimes \id\right) \circ \big[\mathring{\rho} (a_{\eta}^j) + \sum\limits_{\left(\eta\right)} a_{\eta'}^j \otimes a_{\eta''}^k\big] + a_{x_k}^j \otimes a_{\eta}^k$\, , 
\enden 
where $\Delta_{\shuffle}'\left(\eta\right) = \sum\limits_{\left(\eta\right)} \eta' \otimes \eta''$.
\endth

\begpr

Since $\rho\left(V^{+}\right) \subseteq V^{+} \otimes \mathcal{H}$~(c.f.~Lemma~\ref{lem:grade_preservation}), 
\begin{align*}
\mathring{\rho} =  \left(\id \otimes \pi_{V^{+}}\right) \circ \rho
\end{align*}

\begen[(i)]
\item Follows from Proposition~\ref{prop:red_coprod_comp}.
\medskip

\item 
\begin{align*}
\mathring{\rho} \circ \theta_0 &=  \left(\id \otimes \pi_{V^{+}}\right) \circ \rho \circ \theta_0  \\
&= \left(\id \otimes \pi_{V^{+}}\right) \circ \left(\theta_0 \otimes \id\right) \circ \rho \\
&= \left(\theta_0 \otimes \id\right) \circ \left(\id \otimes \pi_{V^{+}}\right) \circ \rho \\
&= \left(\theta_0 \otimes \id\right) \circ \mathring{\rho}.
\end{align*}

\item On $V^{+} \subsetneq \mathcal{B}$,
\begin{align*}
\rho \circ \theta_k\left(a_{\eta}^j\right) 
	= \left(\theta_k \otimes \HAprod\right) \circ \left(\rho \otimes \id_{\mathcal{H}}\right) \circ 
	\Big[\sum_{\left(\eta\right)} a_{\eta'}^j \otimes a_{\eta''}^k + a_{\mathbf{1}}^j \otimes a_{\eta}^k 
	+ a_{\eta}^j \otimes a_{\mathbf{1}}^k\Big],
\end{align*} 
where $\Delta_{\shuffle}'\left(\eta\right) = \sum_{\left(\eta\right)} \eta' \otimes \eta''$ is the reduced unshuffle in $\allpoly$.
\smallskip

Let $\rho\left(a_{\eta}^j\right) = \sum_{\left(\eta\right)}a_{\eta_{(1)}}^j \otimes a_{\eta_{(2)}}^s + a_{\mathbf{1}}^j \otimes a_{\eta}^s$~(for some index family $s$). Thus by Proposition~\ref{prop:red_coprod_comp},

\begin{align*}
\rho \circ \theta_k \left(a_{\eta}^j\right) 
&= \left(\theta_k \otimes \HAprod \right) \circ \Big[\sum_{\left(\eta\right)}\rho\left(a_{\eta'}^j\right) \otimes a_{\eta''}^k + \rho\left(a_{\mathbf{1}}^j\right) \otimes a_{\eta}^k + \rho\left(a_{\eta}^j\right) \otimes a_{\mathbf{1}}^k\Big] \\
&=  \left(\theta_k \otimes \HAprod\right) \circ \Big[\sum_{\left(\eta\right)} a_{\eta'}^j \otimes a_{\mathbf{1}}^s \otimes a_{\eta''}^k + a_{\eta'_{(1)}}^j \otimes a_{\eta''_{(2)}}^s \otimes a_{\eta''}^k + a_{\mathbf{1}}^j \otimes \mathsf{1} \otimes a_{\eta}^k \\
& \quad \quad  \quad \quad + a_{\eta_{(1)}}^j \otimes a_{\eta_{(2)}}^s \otimes a_{\mathbf{1}}^k \Big]\\
&= \left(\theta_k \otimes \id_{\mathcal{B}}\right) \circ \Big[\sum_{\left(\eta\right)} a_{\eta'}^j \otimes a_{\mathbf{1}}^s a_{\eta''}^k + a_{\eta'_{(1)}}^j \otimes a_{\eta''_{(2)}}^sa_{\eta''}^k + a_{\mathbf{1}}^j \otimes a_{\eta}^k + a_{\eta_{(1)}}^j \otimes a_{\eta_{(2)}}^s a_{\mathbf{1}}^k\Big]\\
&= \Big[\left(\theta_k \otimes \id_{\mathcal{B}}\right) \circ \Big[\sum_{\left(\eta\right)} a_{\eta'}^j \otimes a_{\mathbf{1}}^s a_{\eta''}^k 
	+ a_{\eta'_{(1)}}^j \otimes a_{\eta''_{(2)}}^sa_{\eta''}^k + a_{\eta_{(1)}}^j \otimes a_{\eta_{(2)}}^s a_{\mathbf{1}}^k\Big]\Big] + a_{x_k}^j \otimes a_{\eta}^k
\end{align*}

On passing to $\mathcal{H}$ viz. identify all $a_{\mathbf{1}}^j$ with the unit $\mathsf{1}$,
\begin{align*}
\rho \circ \theta_k \left(a_{\eta}^j\right) = \Big[\left(\theta_k \otimes \id_{\mathcal{H}} \right) \circ \Big[\sum_{\left(\eta\right)} a_{\eta'}^j \otimes a_{\eta''}^k + a_{\eta_{(1)}}^j \otimes a_{\eta_{(2)}}^s + \underbrace{a_{\eta'_{(1)}}^j \otimes a_{\eta''_{(2)}}^sa_{\eta''}^k}_{\ker(\pi_{V^{+}})} \Big]\Big] + a_{x_{k}}^j \otimes a_{\eta}^k.
\end{align*}
Therefore,
\begin{align*}
\mathring{\rho} \circ \theta_k \left(a_{\eta}^j\right) 
= \left(\theta_k \otimes \id \right) \circ \Big[\mathring{\rho} \left(a_{\eta}^j\right) 
	+ \sum\limits_{\left(\eta\right)} a_{\eta'}^j \otimes a_{\eta''}^k\Big] + a_{x_k}^j \otimes a_{\eta}^k
\end{align*}
\enden
\endpr

\medskip

The following are some computation of $\mathring{\rho}$ are (for all $i,j,k,\ell = 1,2,\ldots,m$):
\begin{align*}
\mathring{\rho}(a_{x_ix_0}^j) &= a_{x_i}^j \otimes a_{x_0}^i \\
\mathring{\rho}(a_{x_0x_i}^j) &= 0 \\
\mathring{\rho}(a_{x_ix_j}^k) &= a_{x_i}^k \otimes a_{x_j}^i \\
\mathring{\rho}(a_{x_0^2}^i) &= 0 \\
\mathring{\rho}(a_{x_0^3}^i) &= 0 \\
\mathring{\rho}(a_{x_0x_ix_j}^k) &= a_{x_0x_i}^k \otimes a_{x_j}^i \\
\mathring{\rho}(a_{x_kx_ix_j}^{\ell}) &= a_{x_k}^{\ell} \otimes a_{x_ix_j}^{k} + 2a_{x_kx_i}^{\ell} \otimes a_{x_j}^k + a_{x_kx_j}^{\ell} \otimes a_{x_i}^{k}.
\end{align*} 
\medskip

\begre{~For all $n \geq 1$ and $j=1,2,\ldots,m$:\, $a_{x_0^n}^j \in \ker\left(\mathring{\rho}\right)$.}
\medskip


\subsection*{Graded dual of $\mathring{\rho}$}

The graded dual of $V^{+}$ is identified with proper polynomials $\allproperpoly{m} \subsetneq \allseries{m}$; with dual basis $\eta e_j$ where $\eta \in X^{+}$ and $e_j$ for $j = 1,2, \ldots,m$ are standard unit vectors in $\re^m$ such that $a_{\eta}^j\left(\zeta e_k\right) = \delta_{\eta, \zeta} \delta_{j,k}$ for all $\zeta \in X^{+}$. 
\smallskip

The vector space $\allproperpoly{m}$ is equipped with a magmatic product $\triangleleft: \allproperpoly{}^{\otimes 2} \longrightarrow \allproperpoly{}$ such that for $c,d \in \allproperpoly{m}$
\begin{align*}
	(c \triangleleft d)_i\left(\eta\right) 
		= a_{\eta}^i\left(c \triangleleft d\right) = \mathring{\rho}\left(a_{\eta}^i\right)\left(c\otimes d\right), 
\end{align*}
for all $a_{\eta}^i \in V^{+}$. Dualizing Theorem~\ref{th:rho-lin}, the following theorem is obtained:

\begth\label{th:bullet} For all $c,d \in \allproperpoly{m}$ and $j = 1,2,\ldots,m$.
\begen[(i)]
\item $x_0e_j\triangleleft d = 0$ 
\medskip
\item $x_ke_j\triangleleft d = x_kd_ke_j \quad \quad \forall\, k = 1,2,\ldots,m$.
\medskip
\item $x_0c \triangleleft d = x_0\left(c \triangleleft d \right)$.
\medskip
\item $x_kc \triangleleft d = x_k\left(c \triangleleft d \right) + x_k\left(c \shuffle_k \,  d\right) \quad \quad \forall\, k = 1,2, \ldots,m$.
\enden
\endth

\begpr For $d \in \allproperpoly{m}$ and $\eta \in X^{+}, k \neq 0$
\begin{align*}
a_{x_0\eta}^{j}(x_0e_j \triangleleft d) 
	&= \left(\theta_0 \otimes \id\right) \circ \mathring{\rho}\left(a_{\eta}^j\right)\left(x_0e_j \otimes d\right) = 0 \\
a_{x_k\eta}^j\left(x_0e_j \triangleleft d\right) 
	&= \Big[\left(\theta_k \otimes \id \right) \circ \Big[\mathring{\rho} \left(a_{\eta}^j\right) 
	+ \sum\limits_{\left(\eta\right)} a_{\eta'}^j \otimes a_{\eta''}^k\Big] 
	+ a_{x_k}^j \otimes a_{\eta}^k\Big](x_0e_j \otimes d) = 0.
\end{align*}
\smallskip

Similarly, 
\begin{align*}
a_{x_k\eta}^j\left(x_ke_j \triangleleft d\right) 
	&= \Big[\left(\theta_k \otimes \id \right) \circ \Big[\mathring{\rho} \left(a_{\eta}^j\right) 
	+ \sum\limits_{\left(\eta\right)} a_{\eta'}^j \otimes a_{\eta''}^k\Big] + a_{x_k}^j \otimes a_{\eta}^k\Big](x_ke_j \otimes d)\\
&\quad = a_{\eta}^k\left(d\right) = a_{x_k\eta}^j\left(x_kd_ke_j\right) \\
a_{x_0\eta}^j\left(x_ke_j \triangleleft d\right) &= \left(\theta_0 \otimes \id\right) \circ \mathring{\rho}\left(a_{\eta}^j\right)\left(x_ke_j \otimes d\right) = 0 \\
&\quad = a_{x_0\eta}^j(x_kd_ke_j).
\end{align*}

Let $c,d \in \allproperpoly{m}$ and $\eta \in X^{+}$. Using Theorem~\ref{th:rho-lin},
\begin{align*}
a_{x_0\eta}^i\left(c \triangleleft d\right) 
	&= \left(\theta_0 \otimes \id\right) \circ \mathring{\rho}\left(a_{\eta}^i\right)\left(c \otimes d\right)\\
	&= a_{\eta}^i\left(x_0^{-1}\left(c\right) \triangleleft d\right)\Big[\left(\theta_k \otimes \id \right) \circ \Big[\mathring{\rho} \left(a_{\eta}^j\right) + \sum\limits_{\left(\eta\right)} a_{\eta'}^j \otimes a_{\eta''}^k\Big] + a_{x_k}^j \otimes a_{\eta}^k\Big](c\otimes d).
\end{align*}
\smallskip

For $x_k \neq x_0$, 
\begin{align*}
a_{x_k\eta}^j\left(c \triangleleft d\right) 
	&= \Big[\left(\theta_k \otimes \id \right) \circ \Big[\mathring{\rho} \left(a_{\eta}^j\right) 
	+ \sum\limits_{\left(\eta\right)} a_{\eta'}^j \otimes a_{\eta''}^k\Big] + a_{x_k}^j \otimes a_{\eta}^k\Big](c \otimes d) \\
&= a_{\eta}^j \left(x_k^{-1}\left(c\right) \triangleleft d + x_k^{-1}\left(c\right) \shuffle_k\,d\right)
\end{align*}
Hence,
\begin{align*}
a_{x_0\eta}^i\left(x_0c \triangleleft d\right) &= a_{\eta}^i\left(c \triangleleft d\right) = a_{x_0\eta}^i\left(x_0(c \triangleleft d)\right)\\ 
a_{x_k\eta}^i\left(x_0c \triangleleft d\right) &= 0 = a_{x_k\eta}^i\left(x_0\left(c\triangleleft d\right)\right) \\
a_{x_k\eta}^i\left(x_kc \triangleleft d\right) &= a_{\eta}^i\left(c \triangleleft d + \left(c \shuffle_k d\right)\right)  = a_{x_k\eta}^i\left(x_k\left(c \triangleleft d + \left(c \shuffle_k d\right)\right)\right) \\
a_{x_0\eta}^i\left(x_kc \triangleleft d\right) &= 0 = a_{\eta}^i\left(x_k\left(c \triangleleft d + \left(c \shuffle_k d\right)\right)\right).
\end{align*}
Therefore, for all $c,d \in \allproperpoly{m}$,
\begin{align*}
x_0c \triangleleft d &= x_0\left(c \triangleleft d \right) \\
x_kc \triangleleft d &= x_k\left(c \triangleleft d \right) + x_k\left(c \shuffle_k \,  d\right) \quad \quad \forall\, k = 1,2, \ldots,m.
\end{align*}
\endpr
\medskip

\begre{~The product $\triangleleft$ is the graded dual of linearized coaction map. We show in Section~\ref{sec:com-pre-Lie}, that $\left(\allproperpoly{m}, \triangleleft\right)$ is a right pre-Lie algebra with additional structure.}
\medskip


\subsection*{Graded Dual Right pre-Lie Algebra on $\allproperpoly{m}$}

Since $V^{+}$ is a graded right pre-Lie coalgebra, its graded dual $\allproperpoly{m}$ inherits a graded right pre-Lie algebra structure, where pre-Lie product is denoted by $\bullet$ viz
\begin{align*}
	\mathring{\Delta}a_{\eta}^i\left(c\otimes d\right) = \left(c \bullet d\right)_i \left(\eta\right)
\end{align*}

\begth\label{th:pre-Lie-str} $\left(\allproperpoly{m}, \bullet\right)$ is a graded right pre-Lie algebra such that for all $c,d \in \allproperpoly{m}$
\begeq
\label{eq:pre-Lie-prod}
	c \bullet d = \left(c \triangleleft d \right) + \left(c \shuffle d\right),
\endeq

\endth

\begpr This follows from Theorem~\ref{th:bullet}, \ref{eq:copre-Lie-rel} and the fact that
\begin{align*}
\mathring{\Delta}_{\shuffle}a_{\eta}^i\left(c \otimes d\right) = \left(c \shuffle d\right)_j\left(\eta\right).
\end{align*}
\endpr

The derived Lie algebra is $\mathfrak{g}_{\bullet} := \left(\allproperpoly{m},\left[\cdot,\cdot\right]_{\bullet}\right)$.  
\medskip

\begre{~Let $\mathcal{U}\left(\mathfrak{g}_{\bullet}\right)$ denote the universal enveloping algebra of the Lie algebra $\mathfrak{g}_{\bullet}$. Since $\mathcal{H}$ is a graded commutative Hopf algebra, its graded dual, denoted by $\mathring{\mathcal{H}}$, is a cocommutative Hopf algebra. Thus, by Cartier--Milnor--Moore theorem, $\mathring{\mathcal{H}}$ is isomorphic to $\mathcal{U}\left(\mathfrak{g}_{\bullet}\right)$ as graded Hopf algebras.}
\medskip


\section{com-pre-Lie Structure on $\allproperpoly{m}$ associated with a linear Endomorphism}
\label{sec:com-pre-Lie}

The definition of Foissy's com-pre-Lie algebra is given as follows: Let $V$ be a $k$--vector space where $k$ is a field of characteristic $0$.

\smallskip

\begde~\cite{Foissy_15} 
$\left(\mathcal{A},\oslash,\bullet\right)$ is a (right) com-pre-Lie algebra if 
\begen[(i)]
\item $\left(\mathcal{A}, \oslash\right)$ is an associative and commutative algebra.
\medskip
\item $\left(\mathcal{A}, \bullet\right)$ is a right pre-Lie algebra.
\enden
and for all $a,b,c \in \mathcal{A}$
$$
	\left(a \oslash b\right)\bullet c = \left(a \bullet c\right) \oslash b + a \oslash \left(b \bullet c\right).
$$
\endde

\begth
\label{th:com-pre-Lie-2-pre-Lie} \hfill\newline
\begen

\item If $\left(\mathcal{A},\oslash,\bullet\right)$ is right com-pre-Lie, then $\mathcal{A}$ inherits another right pre-Lie product, denoted by $\diamond$ and defined for all $a,b \in \mathcal{A}$ as  $a \diamond b = a \bullet b + a \oslash b$.

\smallskip

\item $\left(\mathcal{A}, \left[\cdot,\cdot\right]_{\diamond}\right)$ is the derived Lie algebra of right pre-Lie algebra $\left(\mathcal{A}, \diamond\right)$ with $\left[a,b\right]_{\diamond} = \left[a,b\right]_{\bullet}$.
\enden
\endth

\begpr
\begen
 \item Observe that $\left(\mathcal{A}, \diamond\right)$ is a Lie--admissible algebra. Indeed, for all $a,b,c \in \mathcal{A}$
\begin{align*}
\left(a \diamond b\right) \diamond c - a \diamond \left(b \diamond c\right) &= \left(a \bullet b + a \oslash b\right) \diamond c - a \diamond \left(b \bullet c + b \diamond c\right) \\
&= \left(a \bullet b \right) \bullet c + \left(a \oslash b\right) \bullet c + \left(a \bullet b\right) \oslash c + \left(a \oslash b\right) \oslash c \\
&\quad - a \bullet \left(b \bullet c\right) - a \bullet \left(b \oslash c\right) - a\oslash \left(b \bullet c\right) - a \oslash\left(b \oslash c\right) \\
&= \left(a\bullet b \right) \bullet c - a \bullet \left(b \bullet c\right) + \left(a\bullet c\right) \oslash b + \left(a \bullet b\right) \oslash c + a \bullet \left(b \oslash c\right).
\end{align*}
The right--hand side is symmetric in $b$ and $c$.
\smallskip

\item Follows from statement $(1)$.
\enden
\endpr
\medskip

\begre{~There are two pre-Lie products, $\diamond$ and $\bullet$, whose derived Lie algebras are identical and with the difference of a commutative product.}
\medskip

Following Foissy \cite{Foissy_15_linear_end}, let $g \in \End{\re X}$, where $\re X$ is the $\re$-span of the alphabet $X$. A com-pre-Lie algebra on proper polynomials, $\allproperpoly{m}$,  associated with such a $g$ (under some assumption) is defined in Theorem~\ref{th:com-pre-Lie-g}, whose definition generalizes the $\triangleleft$ in Theorem~\ref{th:bullet}. By virtue of Theorem~\ref{th:com-pre-Lie-2-pre-Lie}, $\allproperpoly{m}$ inherits another pre-Lie product, which generalizes the pre-Lie product $\bullet$ defined in Theorem~\ref{th:pre-Lie-str}.
\smallskip

 Let $e_j$ for $j =1,2,\ldots,m$ denote the set of standard unit vectors in $\re^m$. 

\begde Define a magmatic product $\triangleleft$ on the vector space $\allproperpoly{m}$ by induction on the degree of polynomials:
\begin{align}\label{eq:pre-Lie-prod-act}
x_ie_j \triangleleft d &= g(x_i)d_ie_j \nonumber\\
x_ic \triangleleft d &= x_i\left(c \triangleleft d\right) + g(x_i)\left(c \shuffle_i d\right) 
\end{align}
where $d_0 := 0$, $x_i \in X$ and $c,d \in \allproperpoly{m}$.
\endde

\begth\label{th:diff-over-sh} For all $c,d,h \in \allproperpoly{m}$:
\begen[(i)]
\item $\left(c \shuffle d\right) \triangleleft h = \left(c \triangleleft h\right) \shuffle d + c \shuffle \left(d \triangleleft h\right)$.
\medskip
\item $\left(c \shuffle_k d\right) \triangleleft h = \left(c \triangleleft h\right) \shuffle_k d + c \shuffle_k \left(d \triangleleft h\right)$.
\enden 
\endth

\begpr
\begen[(i)]
\item  The proof is by induction on $n = \deg\left(c\right) + \deg\left(d\right)$. For base case ($n=2$), let $c = x_i\mathbbm{1}$ and $d = x_j\mathbbm{1}$.
\begin{align*}
\left(x_i\mathbbm{1} \shuffle x_j\mathbbm{1}\right) \triangleleft h
	&= \left(x_ix_j\mathbbm{1} \triangleleft h \right) + \left(x_jx_i\mathbbm{1} \triangleleft h\right) \\
	&= x_i\left(x_j\mathbbm{1} \triangleleft h\right) + g(x_i)\left(x_j\mathbbm{1} \shuffle_i h\right) 
		+ x_j\left(x_i\mathbbm{1} \triangleleft h\right) + g(x_j)\left(x_i\mathbbm{1} \shuffle_i h\right) \\
	&= x_ig(x_j)h_j\mathbbm{1} + g(x_i)\left(x_j\mathbbm{1} \shuffle_i h\right) 
		+ x_jg(x_i)h_i \mathbbm{1} + g(x_j)\left(x_i\mathbbm{1} \shuffle_j h\right) \\
	&= \left(g(x_i)h_i\mathbbm{1} \shuffle x_j \mathbbm{1}\right) + \left(x_i\mathbbm{1} \shuffle g(x_j)h_j\mathbbm{1}\right) \\
	&= \left(x_i\mathbbm{1} \triangleleft h \right) \shuffle x_j\mathbbm{1} 
		+ x_i\mathbbm{1} \shuffle \left(x_j\mathbbm{1} \triangleleft h\right).
\end{align*}
\smallskip
Assume the hypothesis is true up to $n = N$; let $\deg\left(x_ic\right) + \deg\left(x_jd\right) = N+1$; then
\allowdisplaybreaks
\begin{align*}
\left(x_ic \shuffle x_jd\right) \triangleleft h 
	&= \left(x_i\left(c \shuffle x_jd\right) \triangleleft h\right) + \left(x_j\left(x_ic\shuffle d\right)\triangleleft h\right)\\
	&= x_i\left(\left(c \shuffle x_jd\right)\triangleleft h\right) + g(x_i)\left(\left(c \shuffle x_jd\right)\shuffle_i h\right) \\
	&\quad + x_j\left(\left(x_ic \shuffle d\right)\triangleleft h\right) + g(x_j)\left(\left(x_ic \shuffle d\right)\shuffle_j h\right)\\
	&= x_i\left(\left(c \triangleleft h\right) \shuffle x_jd + c \shuffle \left(x_jd \triangleleft h\right)\right) 
		+ x_j\left(\left(x_ic \triangleleft h\right)\shuffle d + x_ic \shuffle \left(d \triangleleft h\right)\right) \\
	&\quad \quad + g(x_i)\left(\left(c \shuffle x_jd\right)\shuffle_i h\right) + g(x_j)\left(\left(x_ic \shuffle d\right) \shuffle_j h\right) \\
&= x_i\left(\left(c \triangleleft h\right)\shuffle x_j d\right) + x_j\left(x_i\left(c \triangleleft h\right) \shuffle d\right) \\
& \quad \quad  + g(x_i)\left(\left(c \shuffle_i h\right) \shuffle x_jd\right) + x_j\left(g(x_i)\left(c \shuffle_i h\right) \shuffle d\right) \\
&\quad x_i\left(c \shuffle x_j\left(d \triangleleft h\right)\right) + x_j\left(x_ic \shuffle \left(d \triangleleft h\right)\right)  \\
&\quad \quad + x_i\left(g(x_j)\left(d \shuffle_j h\right)\right) + g(x_j)\left(x_ic \shuffle \left(d \shuffle_j h\right)\right) \\
&= \left(x_ic \triangleleft h\right) \shuffle x_jd + x_ic \shuffle \left(x_jd \triangleleft h\right). 
\end{align*}

\item Observe that $c \shuffle_k d = c\shuffle d_k\mathbbm{1}$ and thus using $\left(i\right)$,
\begin{align*}
\left(c \shuffle_k d\right) \triangleleft h &= \left(c \shuffle d_k\mathbbm{1}\right) \triangleleft h \\
&= \left(c \triangleleft h \right) \shuffle d_k\mathbbm{1} +  c \shuffle \left(d_k\mathbbm{1} \triangleleft h\right) \\
&= \left(c \triangleleft h \right) \shuffle_k d + c \shuffle \left(d \triangleleft h\right)_k\mathbbm{1}\\
&= \left(c \triangleleft h \right) \shuffle_k d + c \shuffle_k \left(d \triangleleft h\right).
\end{align*}
\enden 
\endpr

\medskip

\begth
\label{th:com-pre-Lie-g} 
$\left(\allproperpoly{m},g,\triangleleft\right)$ is a pre-Lie algebra if and only if $g$ is of the form $g(x_i) = \alpha_i x_i + \beta_i x_0$, for all $i = 1,2,\ldots,m$.
\endth

\begpr It is enough to prove that
\begin{align}\label{eq:identity}
\left(c \triangleleft d\right) \triangleleft h - c \triangleleft \left(d \triangleleft h\right) 
\end{align}
is symmetric in $d$ and $h$ for all $c,d,h \in \allproperpoly{m}$.
\smallskip
Assume that $g(x_i) = \alpha_i x_i + \beta_i x_0$ for $i \neq 0$ and $g(x_0)= \sum_j a_{0j}x_j$. 
\smallskip
The proof is by induction on $n = \deg(c)$. For the base case $(n=1)$; assume $c = x_i\mathbbm{1}$. If $x_i = x_0$ then \rref{eq:identity} is trivial, otherwise
\begin{align*}
\left(x_i\mathbbm{1} \triangleleft d \right) \triangleleft h - x_i\mathbbm{1} \triangleleft \left(d \triangleleft h\right) &= \left(\alpha_i x_id_i\mathbbm{1}\right) \triangleleft h - \alpha_i x_i\left(d \triangleleft h\right)_{i}\mathbbm{1} \\
&= \alpha_i\left(d_i\mathbbm{1} \triangleleft h\right) + \alpha_i x_i\left(d_i \mathbbm{1} \shuffle_i h\right) - \alpha_i x_i\left(d \triangleleft h\right)_i\mathbbm{1} \\
&= \alpha_i x_i \left(d \shuffle h\right)_i \mathbbm{1}
\end{align*}
which is symmetric in $h$ and $d$. Assume that the result is true for $n \leq N$; let $c = x_ic'$ where $\deg(c') = N$. If $x_i = x_0$
\begin{align*}
\left(x_0c' \triangleleft d \right) \triangleleft h - x_0c' \triangleleft \left(d \triangleleft h\right) = x_0\left(c' \triangleleft d\right) \triangleleft h - x_0\left(c' \triangleleft \left(d \triangleleft h\right)\right) = 0.
\end{align*}
For the case $x_i \neq x_0$,
\begin{align*}
\lefteqn{\left(x_ic' \triangleleft d\right) \triangleleft h - x_ic' \triangleleft \left(d \triangleleft h\right)}\\ 
&=  \left(x_i\left(c' \triangleleft d\right) + \alpha_i x_i \left(c' \shuffle_i d\right)\right) \triangleleft h  
		- \left(x_i\left(c'\triangleleft \left(d \triangleleft h\right)\right)  
		- \alpha_ix_i\left(c' \shuffle_i \left(d \triangleleft h\right)\right)\right)\\ 
&= x_i\left(c' \triangleleft d\right)\triangleleft h + \alpha_i x_i\left(c' \shuffle_i d\right)\triangleleft h 
	- x_i \left(c' \triangleleft \left(d \triangleleft h\right)\right) - \alpha x_i\left(c' \shuffle_i \left(d \triangleleft h\right)\right) \\
&= x_i\left(\left(c' \triangleleft d\right) \triangleleft h\right) + \alpha_i x_i\left(\left(c' \triangleleft d\right)\shuffle_i h\right) 
	+ \alpha_i x_i \left(\left(c' \shuffle_i d\right) \triangleleft h \right) + \alpha_i^2x_i\left(\left(c' \shuffle_i d\right)\shuffle_i h\right) \\
&\quad \quad - x_i\left(c' \triangleleft \left(d \triangleleft h\right)\right) - \alpha_i x_i\left(c' \shuffle_i \left(d \triangleleft h\right)\right) \\
&= x_i\left(\left(c'\triangleleft d\right)\triangleleft h - c' \triangleleft \left(d \triangleleft h\right)\right) + \alpha_i^2x_i\left(c' \shuffle_i \left(d \shuffle h\right)\right) \\
&\quad +\alpha_ix_i \left(\left(c' \triangleleft h\right) \shuffle_i d\right) + \alpha_ix_i \left(\left(c' \triangleleft d\right) \shuffle_i h\right) \quad \text{(c.f.~Theorem~ \ref{th:diff-over-sh})}.   
\end{align*}
is symmetric in $d$ and $h$.

Suppose $g(x_i) = a_{i0}x_0 + \sum_{j=1}^{m} a_{ij}x_j$ for $i,j = 1,2,\ldots,m$ such that $\exists\, a_{ij} \neq 0$ for some $i \neq j$.
\begin{align*}
\lefteqn{\left(x_i\mathbbm{1} \triangleleft d \right) \triangleleft h - x_i \mathbbm{1} \triangleleft \left(d \triangleleft h\right)
= \left(g(x_i)d_i\mathbbm{1}\right)\triangleleft h - g(x_i)\left(d \triangleleft h\right)_i\mathbbm{1}} \\
&= g(x_i)\left(d_i\mathbbm{1} \triangleleft h\right) - g(x_i)\left(d_i\mathbbm{1} \triangleleft h\right) + a_{i0}g(x_0)\cancel{\left(d_i\mathbbm \shuffle h_0\right)} 
  +\sum_{j = 1}^m a_{ij}g(x_j)\left(d_i\mathbbm{1} \shuffle_j h\right)  \\
&= \sum_{j=1}^m a_{ij}g(x_j)\left(d_i \shuffle h_j\right)\mathbbm{1} 
\end{align*}
The right-hand is symmetric in $d$ and $h$ only if $i=j$.
\endpr
\medskip

With $X=\{x_0,x_1,x_2,\ldots,x_m\}$ in natural order, the matrix representation of the endomorphism $g$, denoted by $\left[g\right]_{X}$ for which $\left(\allproperpoly{m},g,\triangleleft\right)$ becomes a right pre-Lie algebra is 
\begin{align}\label{eq:matrix-form-g}
\left[g\right]_{X} &= \begin{pmatrix}
a_{00} & a_{01} & a_{02} & a_{03} & \dots & a_{0m} \\
a_{10} & a_{11} & 0 &  0 & \ldots & 0 \\
a_{20} & 0 & a_{22} & 0 & \ldots & 0 \\
\vdots & \vdots  & \vdots & \ddots & \vdots \\
a_{m0} & 0 & 0 & 0 & \ldots & a_{mm}
\end{pmatrix}
\end{align}
\smallskip
where the non--zero (not necessarily zero) elements $a_{ij} \in \re$. The submatrix of $\left[g\right]_{X}$ when restricted to  $X\setminus\{x_0\}$ is a diagonal matrix. 
\medskip

The following theorem is a result of Theorem~\ref{th:diff-over-sh} and Theorem~\ref{th:com-pre-Lie-g}.

\begth\label{th:new-pre-Lie} For $g \in \End{\re X}$ whose matrix representation is of the form in \rref{eq:matrix-form-g},
\begen[(i)]
\item $\left(\allproperpoly{m},g, \shuffle,\triangleleft\right)$ is a right com-pre-Lie algebra. Thus, $\left(\allproperpoly{m}, \bullet\right)$ is a right pre-Lie algebra where $c \bullet d = \left(c \triangleleft d\right) + \left(c \shuffle d\right)$ for all $c,d \in \allproperpoly{m}$.

\item The derived Lie algebras of both right pre-Lie algebras $\left(\allproperpoly{m}, g, \triangleleft\right)$ and $\left(\allproperpoly{m}, \bullet\right)$ are identical. 
\enden
\endth
\medskip


\subsection{pre-Lie in Multiplicative Feedback: A Special Case} \hfill\newline
The product $\triangleleft$ on $\allproperpoly{m}$, defined in Theorem~\ref{th:bullet} obtained as graded dual of linearized coaction map $\mathring{\rho}$ is a special case of $\left(\allproperpoly{m},g, \triangleleft\right)$ with block matrix representation of $g$ as:
\begin{align*}
\left[g\right]_{X} &= \begin{pmatrix}
\begin{matrix}
0
\end{matrix}  & \rvline & 
\begin{matrix}
0 &  \ldots & 0 
\end{matrix} \\
\hline 
\begin{matrix}
0 \\
\vdots \\ 
0
\end{matrix} & \rvline & \LARGE{\mathcal{I}_{m}}
\end{pmatrix}
\end{align*}
where $\mathcal{I}_m$ is the $m \times m$ identity matrix.
\medskip

The following result is a thus a consequence of Theorem~\ref{th:new-pre-Lie}.

\begth\label{th:tri-is-pre-Lie} \hfill
\begen
\item The product $\triangleleft : \allproperpoly{m}^{\otimes 2} \longrightarrow \allproperpoly{m}$, defined in Theorem~\ref{th:bullet} is  (right) pre-Lie. 
\medskip

\item $\left(\allproperpoly{m},\shuffle,\triangleleft\right)$ is a com-pre-Lie algebra.
\enden 
\endth

\begre{~Thus, the Hopf algebra $\mathcal{H}$ can be reconstructed by taking the graded dual of universal enveloping Hopf algebra of the right pre-Lie algebra  $\left(\allproperpoly{m}, \triangleleft\right)$ obtained via Guin--Oudum construction~\cite{GuinOudom2008preLie}.}
\medskip

\subsection*{$m = 1$: A Special Case}
For $m = 1$, then $X = \{x_0,x_1\}$, then $\left(\allproperpoly{},g,  \triangleleft\right)$ is a right com-pre-Lie algebra for any linear endomorphism $g$ as from $\rref{eq:matrix-form-g}$,
\begin{align*}
\left[g\right]_{X} = \begin{pmatrix}
a_{00} & a_{01} \\
a_{10} & a_{11}
\end{pmatrix}
\end{align*}
which are all arbitrary elements $a_{ij} \in \re$. For more details on the com-pre-Lie algebra $\left(\allproperpoly{},g,\triangleleft\right)$, we refer to Foissy~\cite{Foissy_15_linear_end}. As a particular example, for the nilpotent matrix
\begin{align*}
\left[g\right]_{X} = \begin{pmatrix}
0 & 1 \\
0 & 0
\end{pmatrix}
\end{align*}
we obtain the com-pre-Lie (sub)algebra $\left(\allproperpoly{},g,\triangleleft\right)$ corresponding to SISO (single--input, single--output) additive feedback interconnection. It is a com-pre-Lie subalgebra as, the pre-Lie--algebra corresponding to additive feedback is defined on the space of $\allpoly{} \subseteq \allseries{}$ (not just $\allproperpoly{}$). For more details, we refer to Foissy~\cite{Foissy_15}.


\end{document}